# An Improved Mathematical Model of Sepsis: Modeling, Bifurcation Analysis, and Optimal Control Study for Complex Nonlinear Infectious Disease System

Yuyang Chen, Kaiming Bi, Chih-Hang J. Wu, David Ben-Arieh, Ashesh Sinha

## Abstract

Sepsis is a life-threatening medical emergency caused by extreme host immune response to infection, which is a major cause of death worldwide and the second highest cause of mortality in the United States. The immune response is a complicated system. Thus, a more accurate mathematical model is an important tool to study the progression of sepsis. On top of that, researching the optimal control treatment or intervention strategy on the comprehensive sepsis system is key in reducing mortality. For this purpose, first, this paper improves a complex nonlinear sepsis model proposed in our previous work. Then, bifurcation analyses are conducted for each sepsis subsystem to study the model behaviors under some system parameters. The bifurcation analysis results also further indicate the necessity of control treatment and intervention therapy. If the sepsis system is without adding any control under some parameter and initial system value settings, the system will perform persistent inflammation outcomes as time goes by. Therefore, we develop our complex improved nonlinear sepsis model into a sepsis optimal control model, and then use some effective biomarkers recommended in existing clinic practices as optimization objective function to measure the development of sepsis. Besides that, a Bayesian optimization algorithm by combining Recurrent neural network (RNN-BO algorithm) is introduced to predict the optimal control strategy for the studied sepsis optimal control system. The difference between the RNN-BO algorithm from other optimization algorithms is that once given any new initial system value setting (initial value is associated with the initial conditions of patients), the RNN-BO algorithm is capable of quickly predicting a corresponding time-series optimal control based on the historical optimal control data for any new sepsis patient. To demonstrate the effectiveness and efficiency of the RNN-BO algorithm on solving the optimal control solution on the complex nonlinear sepsis system, some numerical simulations are implemented by comparing with other optimization algorithms in this paper.

## 1. Introduction

Sepsis is defined as life-threatening medical emergency caused by the body's extreme systemic immunological response to infection [1]. If there is no any therapeutic treatment, sepsis will further develop into septic shock, organ dysfunction and ultimately result in death. Sepsis is the major causes of death worldwide, with approximately 48.9 million incident sepsis cases in 2017 and estimated 20% of all global deaths [2]. In the early stage of sepsis, source control and antibiotics is normal therapeutic treatment to treat sepsis patients [3]. Some patients are benefit from the early administration of antibiotics [4]. If the patients present persistent inflammation in the later stage of sepsis when bacterial clearance is finished, some studies reported that the anti-TNF-$\alpha$ treatment is an effective therapy [5, 6]. Successful sepsis treatments involve the timing of control therapy and optimal dosing, delayed administration or improper dosage might lead to detrimental outcomes [7]. Thus, providing optimal treatment (involves timing and dosing of administration of control therapy) is the key in reducing the mortality of sepsis and improving patients' quality of care. In the past, attempts to discover the optimal treatments for sepsis have been focused on clinic trails. However, these attempts took much time to manipulate. Also, patients may present different clinical phenotypes if they perform different pathophysiological mechanisms [4], it raises the difficulty to timely provide the effective and appropriate optimal control or intervention treatment through manipulating clinic trails for patients. Therefore, we attempt to address this challenge by the combining use of Bayesian optimization (BO) algorithm and Recurrent neutral network (RNN) applied to a sufficiently complex, nonlinear, mathematical sepsis model.

There are some previous researches on mathematical sepsis model. In 2004, Kumar *et al.* proposed a simplified sepsis mathematical model, this model contains three equations to roughly describe the dynamics between pathogen, early pro-inflammatory and late pro-inflammatory mediators [8]. In 2006, Reynolds *et al.* proposed a sepsis mathematical model to capture scenarios of inflammatory response to infection, this model presents more details of pro-inflammatory and anti-inflammatory mediators [9]. In 2015, We proposed an 18-equation complex sepsis model [10]. This model considers the basic and key components of sepsis progression incorporating innate with adaptive immunities, which studies details immune response among cell, pro-inflammatory cytokines, and anti-inflammatory cytokines. These mathematical models offer insights into complex dynamic immune response. However, these models do not consider the control or intervention treatment as variables into the system, to study the impact of control treatment on sepsis progression and look for the optimal treatment. In order to achieve our original goal, addressing the challenge and studying the method that can timely generate the optimal treatment, in this paper we are therefore developing our previous model into an optimal control model of sepsis.

To construct the sepsis optimal control model, the primary thing is to determine the practical and controllable parameters of the system. Past clinical studies show that appropriate antibiotics therapy in the early hours of sepsis onset can effectively control the pathogen infection, control the pathogen replication/growth rate and decrease absolute mortality [11, 12, 13]. In addition, during the immune response process, the release of the key pro-inflammatory cytokine such as tumor necrosis factor-$\alpha$ (TNF-$\alpha$) is a double-edged sword in sepsis [14], which release rate can positively or negatively influence the outcomes of sepsis progression [10]. Some experimental studies show that anti-TNF-$\alpha$ therapy contributes to control the release rate of TNF-$\alpha$, effective anti-TNF-$\alpha$ therapy can improve the outcome [15, 16]. Therefore, some key parameters such as the growth rate of pathogen and the release rate of TNF-$\alpha$ are controllable in real world, considering their related controls in the researched optimal control model will be more meaningful. Moreover, the model behaviors under some important system parameters is studied via stability and bifurcation analysis.

Besides the controllable parameters, the objective function for the sepsis optimal control model is also needed to be determined. What is a good biomarker that is well suited for the measure of immune response or development of sepsis? Some important immune system components can be used as biomarkers to detect changes and development of sepsis [17]. Those components can be pro-inflammatory cytokines such as TNF-$\alpha$, interleukin-1 beta (IL-1$\beta$), interleukin-6 (IL-6), interleukin-8, and high mobility group box 1 (HMGB-1) [18, 19, 20]. Those pro-inflammatory components are associated to the clearance of pathogen. Some anti-inflammatory cytokines related to the down-regulation of the immune system also can be used as biomarkers, such as interleukin-10 (IL-10), transforming growth factor-$\beta$ (TGF-$\beta$), IL-1 receptor antagonists (IL-1ra) [18, 21]. In addition, TNF-$\alpha$ is a major pro-inflammatory cytokine and IL-10 is a crucial anti-inflammatory cytokine [22]. Thus, in many clinical practices, the ratio between TNF-$\alpha$ and IL-10 served as measure and biomarker to monitor sepsis progression [22, 23, 24]. Besides the inflammatory cytokines, several activation markers of immune cells have been recommended as biomarkers of sepsis, such as neutrophil and monocyte/macrophage immune cells [18, 25]. According to the activation state and functions, monocyte immune cells can develop into monocyte-derived type-1 macrophage (M1 macrophage) and monocyte-derived type 2 macrophage (M2 macrophage) [26, 27]. M1 macrophage can promote the inflammation, M2 macrophage contributes to inhibit inflammation [27]. M1 macrophage is defined as the up-regulation biomarker and M2 macrophage is the down-regulation biomarker of inflammatory [28]. Thus, it is reasonable and convincing to establish the objective function based on some of these biomarkers.

For quickly providing the optimal control treatment strategy of the sepsis optimal control model that optimizes the objective function, the next step we aimed to develop an optimization algorithm by the combining use of BO algorithm and RNN. Solving the optimal control strategy of disease model can be viewed as a nonlinear optimization of control problem with time-series system [29, 30]. BO algorithm has been demonstrated to be an effective algorithm to optimize the optimal control strategy of the complex time-series disease system in our previous works [29, 31]. By combining RNN is because, RNN is great at learning the past data in sequence [32]. Since the initial value of sepsis system parameters and system state variables are associated with the initial conditions of patients, different initial values may lead to different outcomes, and may have different optimal control strategies. If we always use the BO algorithm to solve the optimal control solution when the sepsis patient is associated with new different initial value, the optimization process may take a lot of time, which may miss the best time for treatment. Our main idea is to use BO algorithm to generate the corresponding time-series optimal control strategies for system with different initial values. Consider those different system initial values and corresponding optimal control strategies as known historical data. Then leverage RNN to learn those historical data to catch the relationship between initial value and the optimal control strategy obtained by BO algorithm. Once given a new initial value associated with patient's initial condition, the RNN-BO algorithm can timely and effectively predict the corresponding time-series optimal control strategy for this patient. The most contribution of RNN-BO algorithm is that it learns the historical data and generate a prediction model. Once the RNN-BO prediction model is ready, the RNN-BO algorithm only takes about 2 seconds to predict the optimal control strategy for any new given initial system values. It doesn't take time to do the optimization iterations anymore.

The remainder of this paper is organized as follows. Section 2 formulates the comprehensive optimal control model from some subsystems. Section 3 studies the model behaviors under various parameter settings via stability and bifurcation analysis. Section 4 presents the optimization scheme that will be used for solving the optimal control strategy of sepsis system. Then Section 5 implements the numerical simulation experiments to evaluate the effectiveness of proposed optimization scheme on researched sepsis model. Finally, Section 6 provides the conclusions and discusses our future work.

## 2. Model formulation

Our sepsis optimal control model is developed based on our previous sepsis mathematical model [10]. This model describes the dynamic immune response of liver injury or infection among pathogen, pro-inflammatory cytokines, anti-inflammatory cytokines, and immune cells. We develop this optimal control model by incorporating three subsystems.

### 2.1 Neutrophil immune response subsystem

Macrophage is one of the innate host's first lines of defense against bacterial pathogens [33]. In the initial stage of infection, once the intruding pathogens are detected, the resident immune cells such as tissue macrophages and hepatic macrophage (also known as Kupffer cells or resident liver macrophages) will migrate to the site of pathogens to remove pathogen and resolve infections [9, 10]. Meanwhile, those macrophages release signal to resting phagocytes such as neutrophil immune cells. Resting phagocytes are activated and reach to the infection site to engulf the pathogens. In the meantime, these activated phagocytes release pro-inflammatory cytokines such as TNF- $\alpha$, IL-6, IL-8. The pro-inflammatory cytokines will active and recruit more resting phagocytes to the infection site to clear the pathogen. The activation and recruitment of neutrophil promote the clearance of pathogen. However, the chemical substances such as reactive oxygen species (ROS) released by neutrophil cells is harmful, which will damage host tissue and accelerate the death of apoptotic hepatocytes [34, 35, 36]. We have developed this innate immune response process occurring in the early stage of infection into a mathematical model in the previous works [10]. In this paper, we call it neutrophil immune response subsystem, which consists of the following:

$$\frac{dP}{dt} = k_{pg}P\left(1-\frac{P}{P_\infty}\right) - r_{pmk}\frac{[P^n]}{[P^n+k_{c1}^n]}M_{kf}P^* - r_{pn}\frac{[P^n]}{[P^n+k_{c2}^n]}(N_f+N_b)P^* \tag{1}$$

$$\frac{dM_{kf}}{dt} = k_{mk}M_{kf}\left(1-\frac{M_{kf}}{K_\infty}\right) + k_{mkub}M_{kb} - \frac{[P^n]}{[P^n+k_{c1}^n]}M_{kf}P^* - u_{mk}M_{kf} \tag{2}$$

$$\frac{dM_{kb}}{dt} = \frac{[P^n]}{[P^n+k_{c1}^n]}M_{kf}P^* - k_{mkub}M_{kb} \tag{3}$$

$$\frac{dT}{dt} = \left(\frac{r_{t1max}M_{kb}}{m_{t1}+M_{kb}}\right)M_{kb} + \left(\frac{r_{t2max}N_b}{m_{t2}+N_b}\right)N_b - u_tT \tag{4}$$

$$\frac{dN_R}{dt} = k_{rd}N_R\left(1-\frac{N_R}{N_S}\right) - r_1N_R(T+P)^* - u_{nr}N_R \tag{5}$$

$$\frac{dN_f}{dt} = r_1N_R(T+P)^* + k_{nub}N_b - \frac{[P^n]}{[P^n+k_{c2}^n]}N_fP^* - u_nN_f \tag{6}$$

$$\frac{dN_b}{dt} = \frac{[P^n]}{[P^n+k_{c2}^n]}N_fP^* - k_{nub}N_b \tag{7}$$

$$\frac{dr_1}{dt} = k_{r1}\left(1+\tanh\left(N_f^*\right)\right) - u_{r1}r_1 \tag{8}$$

$$\frac{dD}{dt} = r_{hn}\frac{[D^n]}{[D^n+k_{c3}^n]}N_fD^*\left(1-\frac{D}{A_\infty}\right) - r_{ah}D \tag{9}$$

where $P$, $M_{kf}$, $M_{kb}$, $T$, $N_R$, $N_f$, $N_b$, $r_1$, $D$ are $t_f$-dimensional system state variables, $t \in [t_1, t_f]$, $t_1$ is the start time and $t_f$ is the end time. They represent the levels of pathogen, free Kupffer cell that is waiting for binding with pathogen, binded Kupffer cell that is binding with pathogen, TNF- $\alpha$, resting neutrophil that is waiting for activation, free activated neutrophil that is activated and is waiting for binding with pathogen, binded activated neutrophil that is binding with pathogen, the rate of resting neutrophil activated under infection, and damaged tissue or dead hepatocytes, respectively. In Eq. (1), $P^*$ represents the pathogen concentration defined as $P^* = \frac{P}{P_\infty}$. In Eq. (8), $N_f^*$ represents the free activated neutrophil concentration as $N_f^* = \frac{N_f}{N_S}$. In Eq. (9), $D^*$represents the damage tissue concentration as $D^* = \frac{D}{A_\infty}$. The rest of symbols are system parameters, their definition and corresponding values for later simulation experiments are summarized in Table 1 shown in Appendix. We refer readers to our previous work [10] to get more details about the construction of this neutrophil immune response subsystem.

### 2.2 Monocyte immune response subsystem

In our previous work [10], we have also constructed the monocyte immune response subsystem. However, the previous work did not consider the further development of monocytes. To better describe the dynamics of immune response, we attempt to improve the monocyte immune response model in this paper.

During the innate immune response process, besides the presence of Kupffer cell and neutrophil phagocyte contributing to the clearance of pathogen, recent works from the literature have already shown that monocyte immune cell is also a key phagocyte [37]. Monocyte is activated and recruited by HMGB-1 and TNF- $\alpha$, which is capable of clearing the pathogen and phagocytizing the aging binded activated neutrophils, it has significant impact on liver inflammation [10, 38, 39, 40]. On the other hand, according to existing literature, HMGB-1 can be released by activated monocytes and necrotic cells (means dead cells in this paper) [40, 41, 42]. Besides the release of HMGB-1, monocytes also release the anti-inflammatory cytokines such as IL-10 [43]. IL-10 contributes to prevent the subsequent tissue damage by inhibiting the activation of phagocytes such as neutrophils and monocytes [44].

About the monocyte development, many experimental evidence indicates that monocytes will develop into monocyte-derived type 1 macrophage (M1 macrophage) when they encounter pathogen, TNF-α, or GM-CSF, then M1 macrophage contributes to kill the pathogens through phagocytosis [45, 46]. During this process, M1 macrophages will release pro-inflammatory cytokines such as TNF- α, IL-6, and IL-12 [47, 48]. Thus, M1 macrophages are inflammatory microphages that can promote inflammation and cause damage to host tissues [47]. In addition, monocytes also will develop into monocyte-type 2 macrophage (M2 macrophage) when they encounter apoptotic T cells, IL-10, or TGF-$\beta$ [45, 46, 47]. M2 macrophages will release anti-inflammatory cytokines IL-10 and TGF- $\beta$ when they phagocytize apoptotic T cells [49]. Thus, M2 macrophages are healing macrophages that plays an important impact on the healing and tissue repair [47]. A simplified mechanism of monocyte development is drawn as Fig. 1.

Due to the immune response of M2 macrophages is associated with T cell, it belongs to adaptive immunity. Therefore, in this monocyte immune response subsystem, we will only consider M1 macrophages. The mathematical expression of M2 macrophage will be constructed in immune system with adaptive immunity shown in Section 2.3. Base on the original model proposed in previous work [10] and the development of monocytes, the monocyte immune response subsystem is revised by remodeling the expression of monocytes as following:

$$\frac{dP}{dt} = k_{pg}P(1-P^*) - r_{pmk}\frac{[P^n]}{[P^n+k_{c1}^n]}M_{kf}P^* - r_{pn}\frac{[P^n]}{[P^n+k_{c2}^n]}\left(N_f+N_b\right)P^* - E_1 \tag{10}$$

$$\frac{dN_f}{dt} = \frac{r_1 N_R (T+P)^*}{(1+\frac{C_A}{C_\infty})} + k_{nub}N_b - \frac{[P^n]}{[P^n+k_{c2}^n]}N_f P^* - u_n N_f \tag{11}$$

$$\frac{dN_b}{dt} = \frac{[P^n]}{[P^n+k_{c2}^n]}N_f P^* - u_{mn}N_b M_f^* - k_{nub}N_b \tag{12}$$

$$\frac{dM_R}{dt} = k_{mr}M_R\left(1-\frac{M_R}{M_S}\right) - r_2 M_R (H+T)^* - u_{mr}M_R \tag{13}$$

$$\frac{dM_f}{dt} = \frac{r_2 M_R (H+T)^*}{(1+\frac{C_A}{C_\infty})} + k_{umb} M_b - E_1 - u_m M_f \tag{14}$$

$$\frac{dM_b}{dt} = E_1 - k_{umb} M_b \tag{15}$$

$$E_1 = \frac{dM_1}{dt} = r_{pm} \frac{[P^n]}{[P^n + k_{c4}^n]} M_f P^* \tag{16}$$

$$\frac{dH}{dt} = \left(\frac{r_{h1max}(M_b+D)}{mh_1+M_b+D}\right)(M_b + D) - u_h H \tag{17}$$

$$\frac{dC_A}{dt} = \left(\frac{r_{camax} M_b}{C_{Ah}+M_b}\right) M_b - u_{ca} C_A \tag{18}$$

where $M_R$, $M_f$, $M_b$, $M_1$, $H$, $C_A$ are $t_f$-dimensional system state variables, $t \in [t_1, t_f]$, $t_1$ is the start time and $t_f$ is the end time. They represent the levels of resting monocyte that is waiting for activation, free activated monocyte that is activated and is waiting for phagocytizing, binded activated monocyte that is involving in the immune response with pathogen and T cells, monocyte-derived type 1 macrophage, HMGB-1, and IL-10, respectively.

Eq. (10) is developed from Eq. (1) by incorporating the clearance effect of monocytes. Eq. (11) is developed from Eq. (6) due to the inhibition of IL-10. Eq. (12) is developed from Eq. (7) by incorporating the phagocytosis effect of monocytes. In Eq. (12), $M_f^*$ represents the free activated monocyte concentration as $M_f^* = \frac{M_f}{M_S}$. Eq. (16) represents the changing number of M1 macrophages due to M1 macrophages phagocytize pathogen, this term is associated with the solid line numbered with ① in Fig. 1. The rest of symbols without mentioned before are system parameters, their definition and corresponding values for later simulation experiments are summarized in Table 1 shown in Appendix.

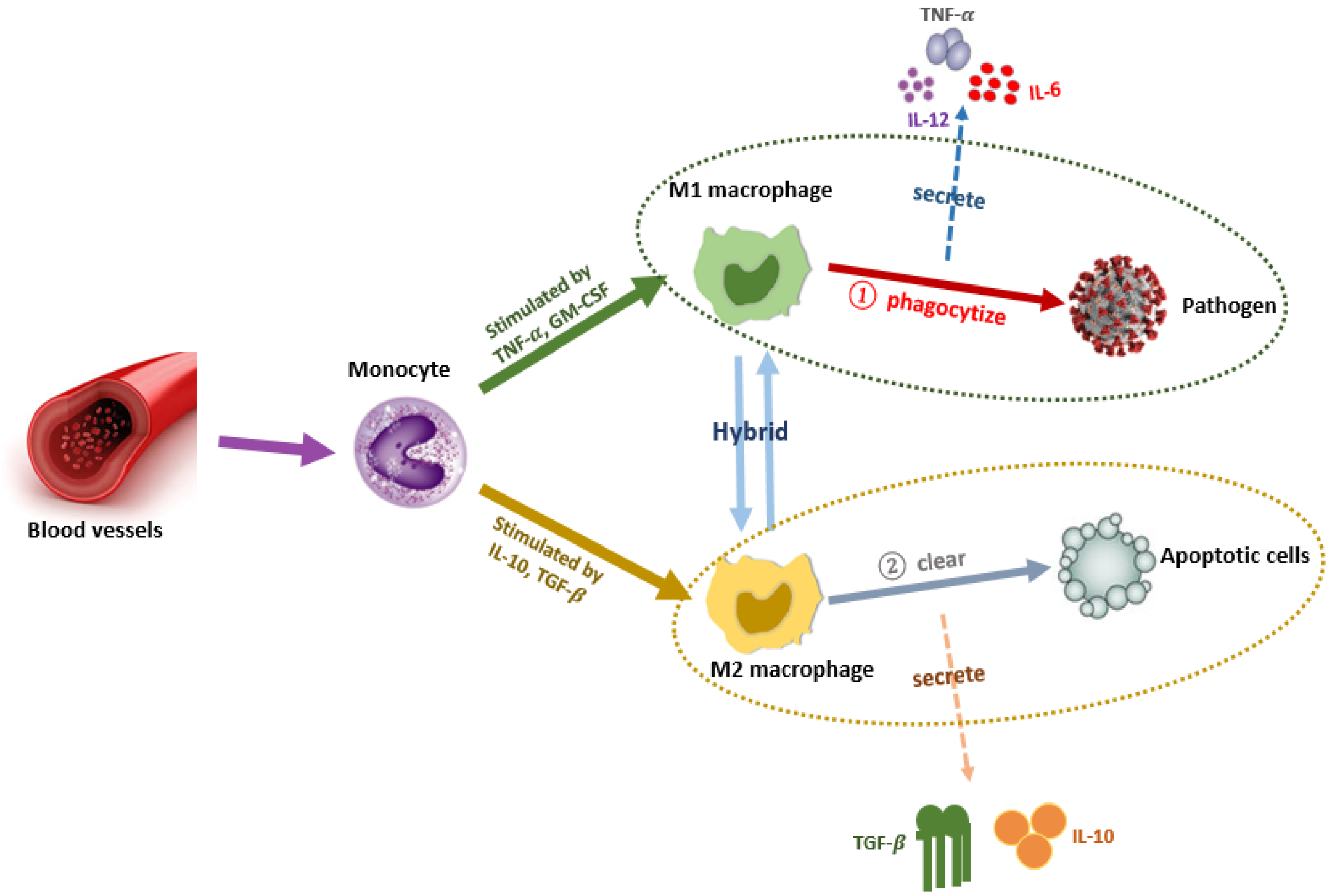

**Fig. 1. Simplified mechanism of monocyte development.**

## 2.3 Immune response system incorporated with adaptive immunity

Innate immunity plays an important role in the clearance of pathogen in the early stage of inflammation. Compared to innate immunity, adaptive immunity is activated in the late stage of inflammation [50]. The dynamics of adaptive immunity is more complicated than innate immunity. To simplify adaptive immunity, in this paper we will remain the model including B cells, and antibodies proposed in [10]. On this basis, this paper will provide the expression of M2 macrophages, and remodel the expression of monocytes. At the same time, this paper will improve the expression of T cells due to some T cell's functions.

The T cells we will study and model in this paper are CD4+ T cell and CD8+ T cell. CD4+ T cells play an important role on clearing the pathogen and achieving a regulated effective immune response to infection [51]. Activated monocytes that phagocytize pathogen is one type of antigen-presenting cells (APCs) [52]. CD4+ T cells are activated and recruited by APCs, APCs also can enhance and recruit more CD4+ T cells [49, 50]. CD4+ T cells that undergo apoptotic are phagocytized by M2 macrophages [53]. The activation of CD8+ T cells go through a major histocompatibility complex class I peptide (MHCI)-TCR mechanism, which is similar to the activation process of CD4+ T cells [10]. CD8+ T cells that undergo apoptotic are phagocytized by M2 macrophages as well [53]. Unlike CD4+ T cells, CD8+ T cells are cytotoxic cells, their primary function is to kill the infected target cells [49, 54]. In previous work, we have modeled the clearance function of CD4+ T cells on pathogen expression and cytotoxic function of CD8+ T cells through the decrease on expressions of binded Kupffer cells, binded activated neutrophils, and binded activated monocytes. However, the previous work doesn't model the clearance function on expression of CD4+ T cells and cytoxic function on expression of CD8+ T cells. Thus, we not only remain the modeling of those functions on pathogen, binded Kupffer cells, binded activated neutrophils, and binded activated monocytes, but also revise the expression of CD4+ T cells and CD8 T cells in this paper. A simplified mechanism of T cells in this paper is drawn as Fig. 2.

Some experimental studies shown have shown that CD4+ T cells are activated by APCs to proliferate and differentiate into $T_H1$ and $T_H2$ effector cells [49, 55]. $T_H1$ and $T_H2$ effector cells can activate B cells to secrete antibodies [50]. The antibodies released by B cells play an important role on the clearance of pathogen at the later stage of inflammation [50, 52].

Base on the original model proposed in previous work [10] and our improvement on monocytes and T cells, the improved immune response system incorporated with adaptive immunity is revised as following:

$$\frac{dP}{dt} = k_{pg}P(1-P^*) - r_{pmk}\frac{[P^n]}{[P^n+k_{c1}^n]}M_{kf}P^* - r_{pn}\frac{[P^n]}{[P^n+k_{c2}^n]}(N_f+N_b)P^* - M_1 - r_{pcd4}\frac{[P^n]}{[P^n+k_{c6}^n]}T_{CD4}P^* - r_{pAb}\frac{[P^n]}{[P^n+k_{c5}^n]}AP^* \tag{19}$$

$$\frac{dM_{kf}}{dt} = k_{mk}M_{kf}\left(1-\frac{M_{kf}}{K_\infty}\right) + k_{mkub}M_{kb} - \frac{[P^n]}{[P^n+k_{c1}^n]}M_{kf}P^* - u_{mk}M_{kf} \tag{20}$$

$$\frac{dM_{kb}}{dt} = \frac{[P^n]}{[P^n+k_{c1}^n]}M_{kf}P^* - r_{Mkbcd8}\frac{[M_{kb}^n]}{[M_{kb}^n+k_{c6}^n]}T_{CD8}M_{kb}^* - k_{mkub}M_{kb} \tag{21}$$

$$\frac{dT}{dt} = \left(\frac{r_{t1max}M_{kb}}{m_{t1}+M_{kb}}\right)M_{kb} + \left(\frac{r_{t2max}N_b}{m_{t2}+N_b}\right)N_b - u_tT \tag{22}$$

$$\frac{dN_R}{dt} = k_{rd}N_R\left(1-\frac{N_R}{N_S}\right) - r_1N_R(T+P)^*/(1+\frac{C_A}{C_\infty}) - u_{nr}N_R \tag{23}$$

$$\frac{dN_f}{dt} = r_1N_R(T+P)^*/(1+\frac{C_A}{C_\infty}) + k_{nub}N_b - \frac{[P^n]}{[P^n+k_{c2}^n]}N_fP^* - u_nN_f \tag{24}$$

$$\frac{dN_b}{dt} = \frac{[P^n]}{[P^n+k_{c2}^n]}N_fP^* - u_{mn}N_bM_f^* - r_{Nbcd8}\frac{[N_b^n]}{[N_b^n+k_{c7}^n]}T_{CD8}N_b^* - k_{nub}N_b \tag{25}$$

$$\frac{dr_1}{dt} = k_{r1}\left(1+\tanh\left(N_f^*\right)\right) - u_{r1}r_1 \tag{26}$$

$$\frac{dD}{dt} = r_{hn}\frac{[D^n]}{[D^n+k_{c3}^n]}N_fD^*(1-\frac{D}{A_\infty}) - r_{ah}D \tag{27}$$

$$\frac{dM_R}{dt} = k_{rm}M_R\left(1 - \frac{M_R}{M_S}\right) - r_2 M_R(H + T + T_{CD4} + T_{CD8})^*/(1 + \frac{C_A}{C_\infty}) - u_{mr}M_R \quad (28)$$

$$\frac{dM_f}{dt} = \frac{r_2 M_R(H+T+T_{CD4}+T_{CD8})^*}{1+\frac{C_A}{C_\infty}} + k_{umb}M_b - E_1 - E_2 - r_{cd4Mb}\frac{[M_f^n]}{[M_f^n + k_{c8}^n]}T_{CD4}M_f^*$$

$$-r_{cd8Mb}\frac{[M_f^n]}{[M_f^n + k_{c8}^n]}T_{CD8}M_f^* - u_m M_f \quad (29)$$

$$\frac{dM_b}{dt} = E_1 + E_2 + r_{cd4Mb}\frac{[M_f^n]}{[M_f^n + k_{c8}^n]}T_{CD4}M_f^* + r_{cd8Mb}\frac{[M_f^n]}{[M_f^n + k_{c8}^n]}T_{CD8}M_f^* - r_{Mbcd8}\frac{[M_b^n]}{[M_b^n + k_{c7}^n]}T_{CD8}M_b^* - k_{umb}M_b \quad (30)$$

$$E_1 = \frac{dM_1}{dt} = r_{pm}\frac{[P^n]}{[P^n + k_{c4}^n]}M_f P^* \quad (31)$$

$$E_2 = \frac{dM_2}{dt} = k_{cd4M}\frac{[T_{CD4}^n]}{[T_{CD4}^n + k_{c10}^n]}M_f T_{CD4}^* + k_{cd8M}\frac{[T_{CD8}^n]}{[T_{CD8}^n + k_{c10}^n]}M_f T_{CD8}^* \quad (32)$$

$$\frac{dH}{dt} = \left(\frac{r_{h1max}(M_b + D)}{mh_1 + M_b + D}\right)(M_b + D) - u_h H \quad (33)$$

$$\frac{dC_A}{dt} = \left(\frac{r_{camax}M_b}{C_{Ah} + M_b}\right)M_b - u_{ca}C_A \quad (34)$$

$$\frac{dT_{CD4}}{dt} = k_{cd4}T_{CD4}\left(1 - \frac{T_{CD4}}{T_{CD4\infty}}\right) + r_{cd4Mb}\frac{[M_f^n]}{[M_f^n + k_{c8}^n]}T_{CD4}M_f^* - k_{cd4M}\frac{[T_{CD4}^n]}{[T_{CD4}^n + k_{c10}^n]}M_f T_{CD4}^*$$

$$-r_{pcd4}\frac{[P^n]}{[P^n + k_{c6}^n]}T_{CD4}P^* - u_{cd4}T_{CD4} \quad (35)$$

$$\frac{dT_{CD8}}{dt} = k_{cd8}T_{CD8}\left(1 - \frac{T_{CD8}}{T_{CD8\infty}}\right) + r_{cd8Mb}\frac{[M_f^n]}{[M_f^n + k_{c8}^n]}T_{CD8}M_f^* - k_{cd8M}\frac{[T_{CD8}^n]}{[T_{CD8}^n + k_{c10}^n]}M_f T_{CD8}^*$$

$$-r_{Mkbcd8}\frac{[M_{kb}^n]}{[M_{kb}^n + k_{c6}^n]}T_{CD8}M_{kb}^* - r_{Nbcd8}\frac{[N_b^n]}{[N_b^n + k_{c7}^n]}T_{CD8}N_b^* - r_{Mbcd8}\frac{[M_b^n]}{[M_b^n + k_{c7}^n]}T_{CD8}M_b^* - u_{cd8}T_{CD8} \quad (36)$$

$$\frac{dB}{dt} = k_B B\left(1 - \frac{B}{B_\infty}\right) + r_{Bt}\frac{[B^n]}{[B^n + k_{c9}^n]}T_{CD4}B^* - u_B B \quad (37)$$

$$\frac{dA}{dt} = \left(\frac{r_{Abmax}B}{m_{Ab} + B}\right)B - u_{Ab}A \quad (38)$$

where $M_2$, $T_{CD4}$, $T_{CD8}$, $B$, $A$ are $t_f$-dimensional system state variables, $t \in [t_1, t_f]$, $t_1$ is the start time and $t_f$ is the end time. They represent the levels of monocyte-derived type 2 macrophage, CD4+ T cell, CD8+ T cell, B cell, and Antibodies, respectively. In Eq. (21), $M_{kb}^*$ represents the binded Kupffer cell concentration defined as $M_{kb}^* = \frac{M_{kb}}{K_\infty}$. In Eq. (8), $N_b^*$ represents the binded activated neutrophil concentration as $N_b^* = \frac{N_b}{N_S}$. In Eq. (25), $M_b^*$ represents the binded activated monocytes concentration as $M_b^* = \frac{M_b}{M_S}$. In Eq. (32), $T_{CD4}^*$ represents the CD4+ T cell concentration as $T_{CD4}^* = \frac{T_{CD4}}{T_{CD4\infty}}$, $T_{CD8}^*$ represents the CD8+ T cell concentration as $T_{CD8}^* = \frac{T_{CD8}}{T_{CD8\infty}}$. In Eq. (38), $B^*$ represents B cell concentration as $B^* = \frac{B}{B_\infty}$. Eq. (32) represents the changing number of M2 macrophages due to M2 macrophages phagocytize apoptotic T cells, this term is associated with the solid line numbered with ② in Fig. 1.

In Eq. (35), the first term represents the recruiting process of CD4+ T cells during adaptive immunity, which is associated with the solid line numbered with ① in Fig. 2. The second term represents the increasing number of CD4+

T cells that are enhanced by APCs, which is associated with the solid line numbered with ③ in Fig. 2. The third term represents the decreasing number of CD4+ T cells since the apoptotic CD4+ T cells are phagocytized by monocytes, which is associated with the solid line numbered with ⑦ in Fig. 2. The fourth term represents the decreasing number of CD4+ T cells since they are binding with pathogen and kill pathogen, which is associated with the solid line numbered with ⑤ in Fig. 2. The fifth term represents the decreasing number of CD4+ T cells due to normal degradation, which is associated with the solid line numbered with ⑥ in Fig. 2.

In Eq. (36), the first term represents the recruiting process of CD8+ T cells during adaptive immunity, which is associated with the solid line numbered with ② in Fig. 2. The second term represents the increasing number of CD8+ T cells that are enhanced by APCs, which is associated with the solid line numbered with ④ in Fig. 2. The third term represents the decreasing number of CD8+ T cells since the apoptotic CD8+ T cells are phagocytized by monocytes, which is associated with the solid line numbered with ⑧ in Fig. 2. The fourth, fifth, sixth terms represent the decreasing number of CD8+ T cells since CD8+ T cells are binding with Kupffer cells, neutrophils, monocytes and kill them, which are associated with the solid lines numbered with ⑨, ⑩, ⑪ in Fig. 2, respectively. The seventh term represents the decreasing number of CD8+ T cells due to normal degradation, which is associated with the solid line numbered with ⑫ in Fig. 2.

The rest of symbols without mentioned before are system parameters, their definition and corresponding values for later simulation experiments are summarized in Table 1 shown in Appendix.

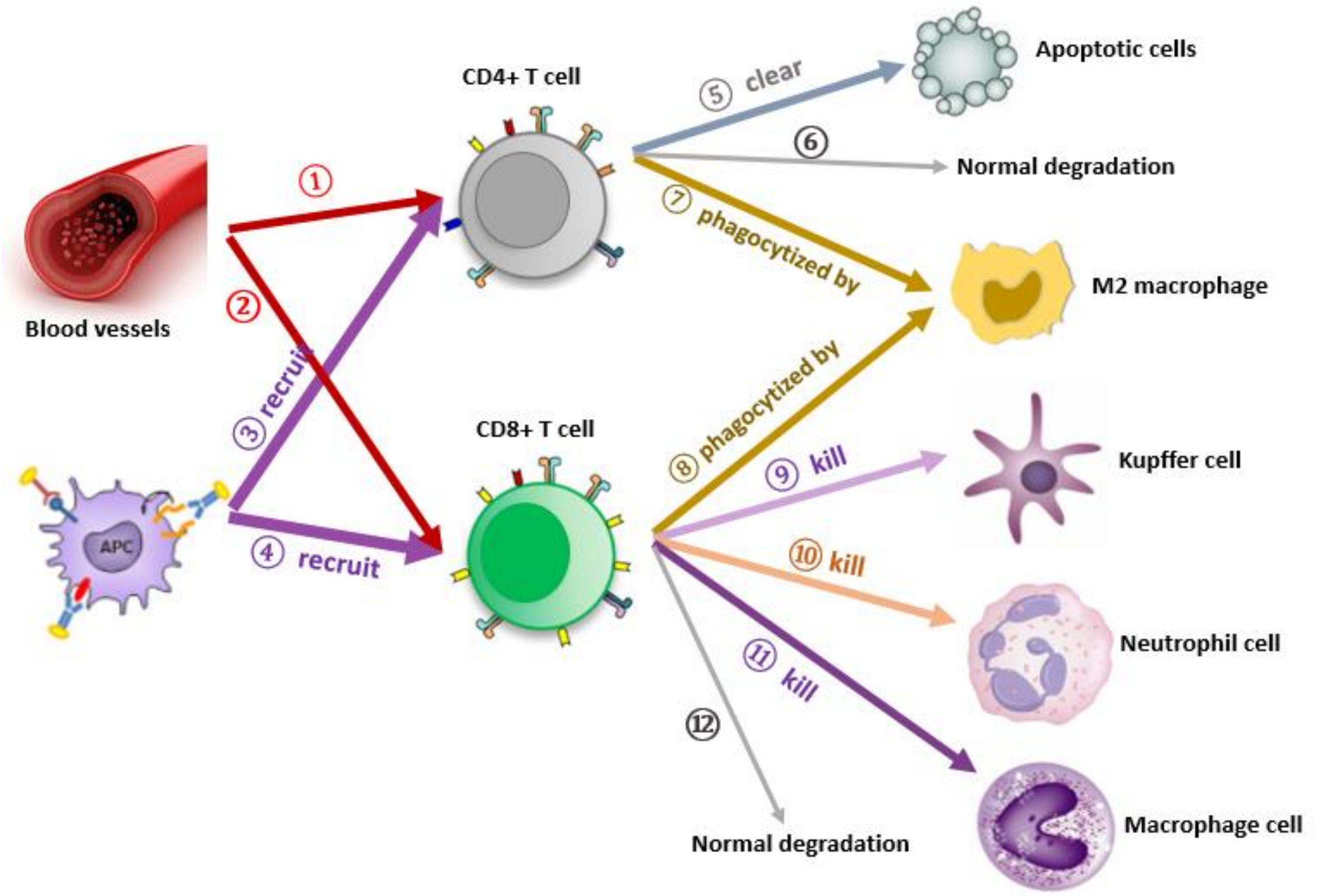

**Fig. 2. Simplified mechanism of T cells.**

## 3. Bifurcation analysis

To study the model dynamics behaviors under various parameter settings, we will conduct the bifurcation analysis for each subsystem in this section. Bifurcation is the qualitative behavior change (change in number or numerical value of equilibrium points) of the system by varying parameters [56]. The objective of bifurcation analysis is to study and identify the key parameters in sepsis development. In this paper we will use numerical analysis to realize bifurcation analysis due to the complexity of sepsis system. Since our current nonlinear sepsis model is too complicated, there is no existing programming tools or packages that can directly solve the bifurcation diagrams of system. Bifurcation value is a value of the equilibrium point moving from stable equilibrium to unstable equilibrium [57]. Therefore, we will start the bifurcation analysis by varying the values of key parameters, then plot all equilibrium points over the key parameters. The bifurcation will be intuitively and clearly caught. In this paper, all bifurcation diagrams are numerically generated by Python.

### 3.1 Bifurcation analysis in neutrophil subsystem

The parameters we analyze in neutrophil subsystem are $k_{pg}$, $r_{pn}$, and $u_n$. For each parameter, only the system state variables with obvious equilibrium behavior are presented. The bifurcation diagrams of neutrophil subsystem are shown in Fig. 3.

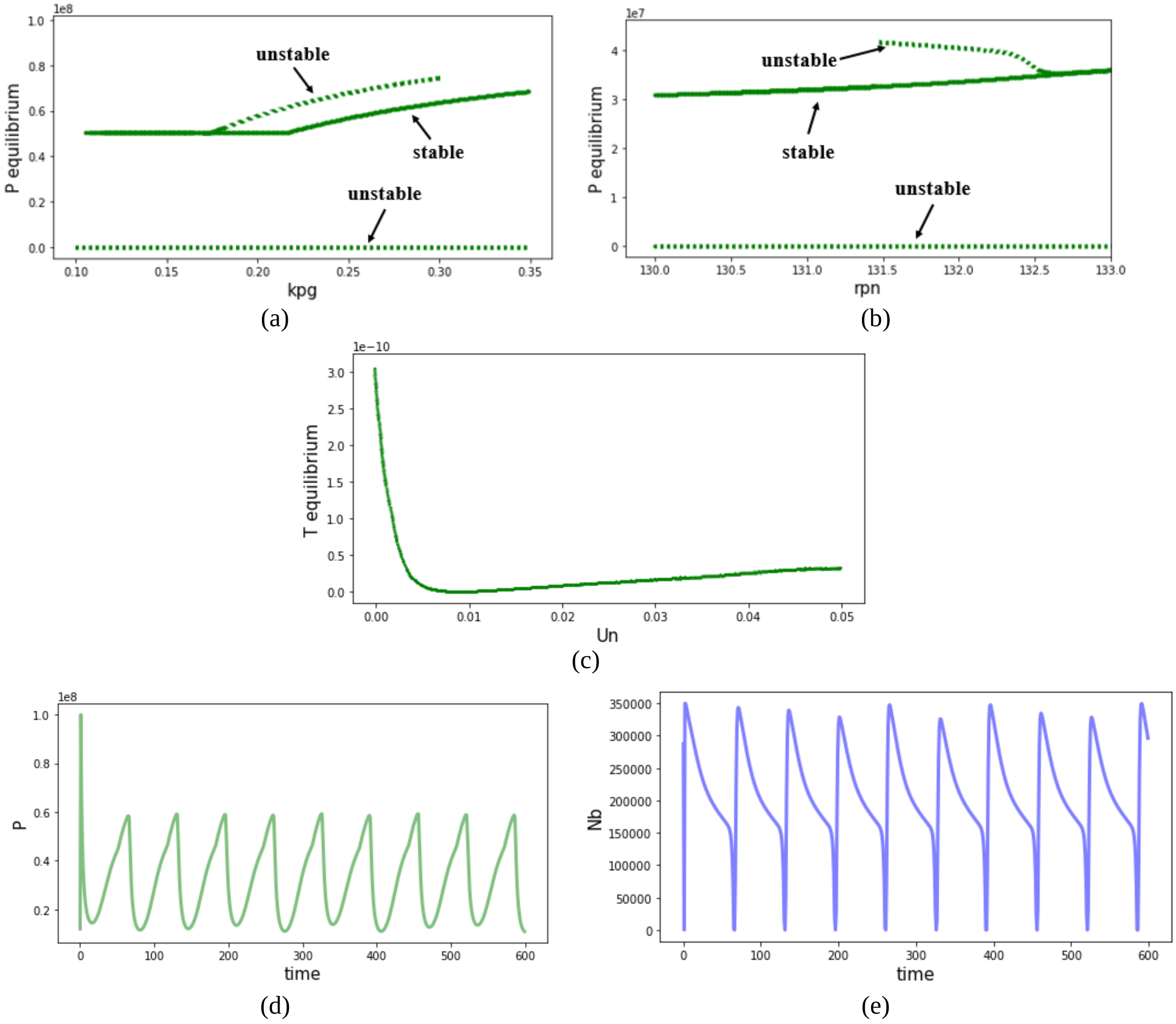

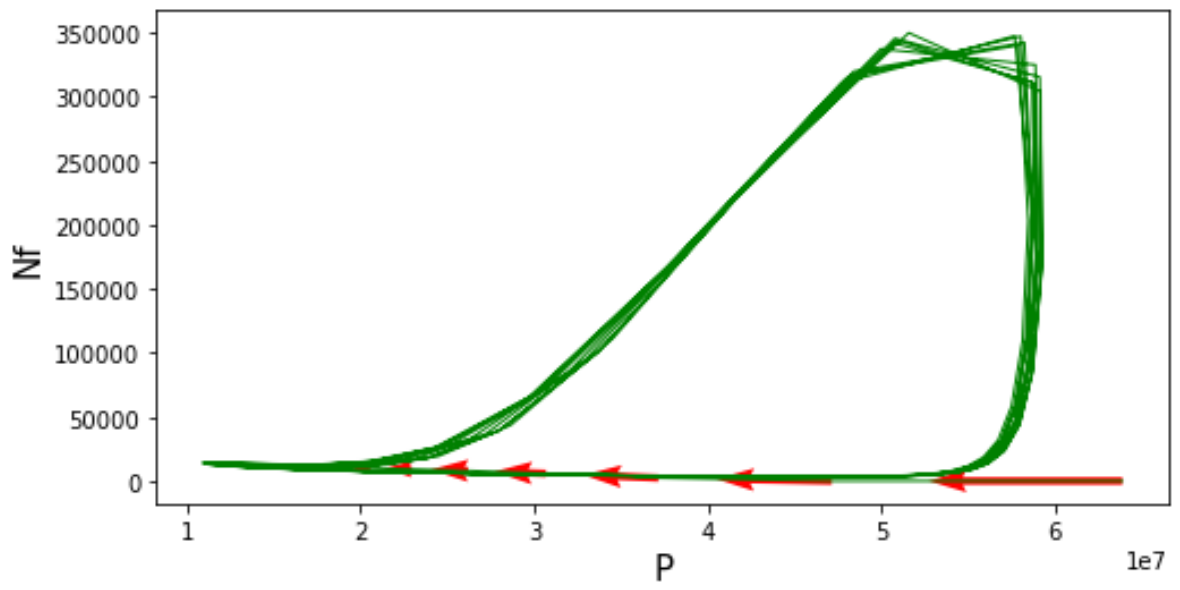

(f)

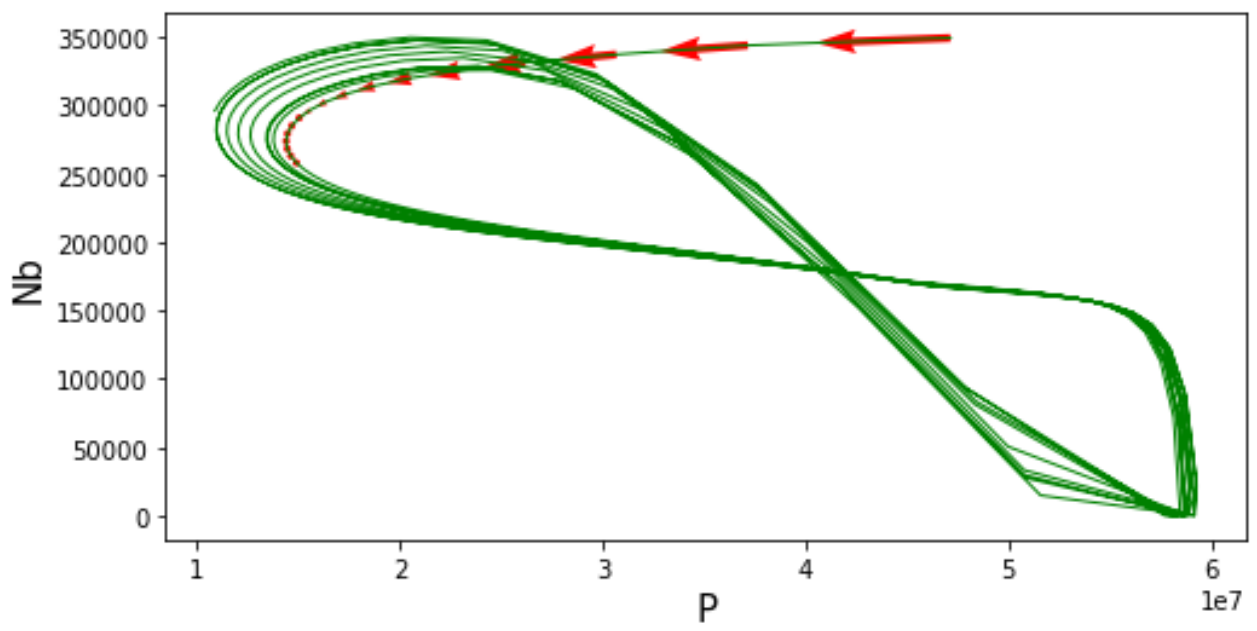

(g)

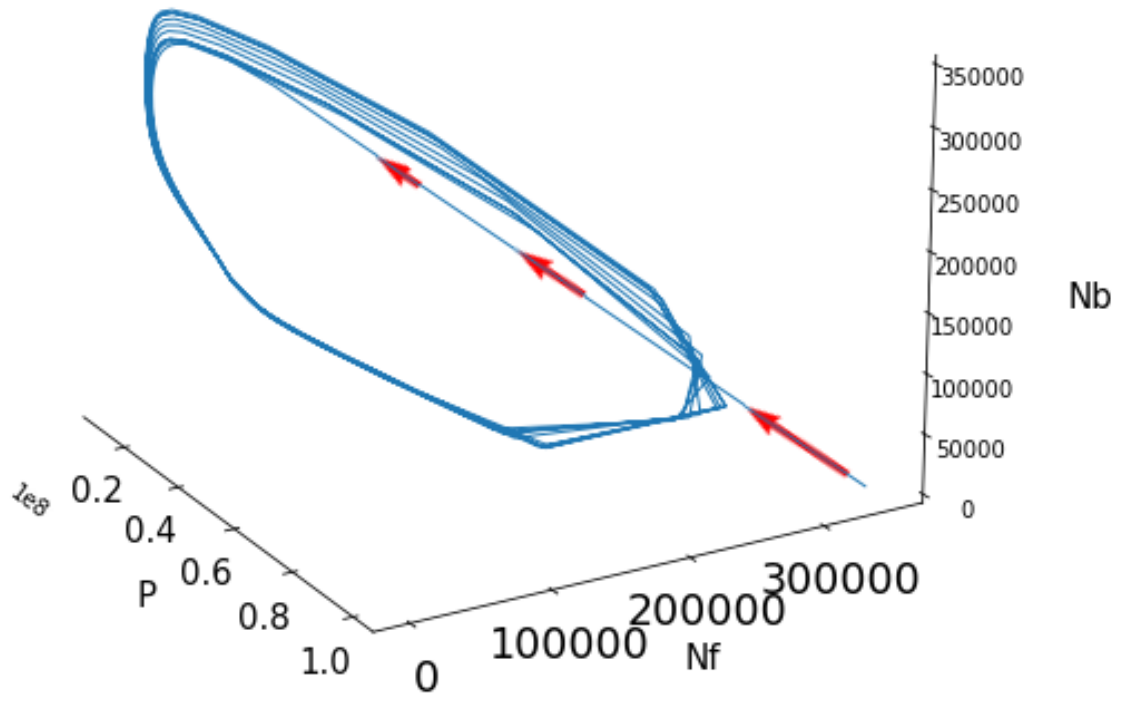

(h)

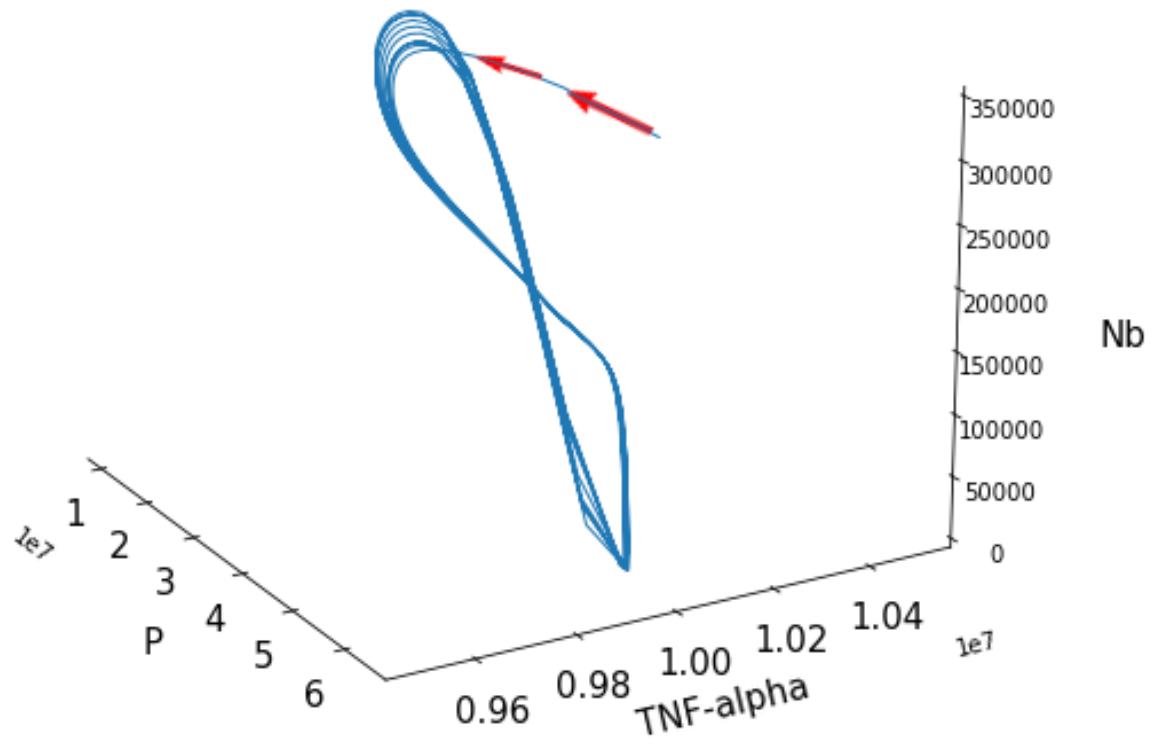

(i)

**Fig. 3. (a) Numerical equilibrium curve of pathogen related to parameter $k_{pg}$ in neutrophil subsystem.**
**(b) Numerical equilibrium curve of pathogen related to parameter $r_{pn}$ in neutrophil subsystem.**
**(c) Numerical equilibrium curve of TNF-$\alpha$ related to parameter $u_n$ in neutrophil subsystem.**
**(d) Oscillation behavior of pathogen in neutrophil subsystem when $k_{pg}$ is equal to 0.1.**
**(e) Oscillation behavior of Nb in neutrophil subsystem when $k_{pg}$ is equal to 0.1.**
**(f) Phase trajectory in P-Nf plane in neutrophil subsystem.**
**(g) Phase trajectory in P-Nb plane in neutrophil subsystem.**
**(h) Phase trajectory in (P, Nf, Nb) space in neutrophil subsystem.**
**(i) Phase trajectory in (P, TNF-$\alpha$, Nb) space in neutrophil subsystem.**

In Fig. 3, X axis represents the parameter values, Y axis represents the equilibrium values of the system state variable. According to the definition of bifurcation, Fig. 3 (a) and (b) both show the changes in the number of equilibrium and the change in the numerical values of equilibrium when the parameter value is change. In (a) and (b), the solid line represents the stable equilibrium, dash line represents the unstable equilibrium. Stable equilibrium means that the points nearing this equilibrium (on both sides of this equilibrium) converge to this equilibrium, unstable equilibrium means that there exist points nearing this equilibrium (on both sides of this equilibrium) diverge from this equilibrium [58]. In Fig. 3 (a), stable equilibrium points of pathogen are observed when system parameter $k_{pg}$ increases from 0.11 to 0.35. At the same range of parameter $k_{pg}$, unstable equilibrium points of pathogen are observed as well. When $k_{pg} = 0.175$, a bifurcation point is identified and new unstable equilibrium point of pathogen are generated as $k_{pg}$ increases from 0.175 to 0.3. In Fig. 3 (b), stable equilibrium points of pathogen are observed when system parameter $r_{pn}$ increases from 130 to 133. At the same range of parameter $r_{pn}$, unstable equilibrium points of pathogen are observed as well. When $r_{pn} = 132.6$, a bifurcation point is identified and new unstable equilibrium point of pathogen are generated as $r_{pn}$ decreases from 132.6 to 131.5. In Fig. 3 (c), the changes on numerical value of equilibrium points of TNF-$\alpha$ is observed by varying the system parameters $u_n$.

Fig. 3 (d) and (e) show the oscillation behaviors of pathogen and $N_b$ when $k_{pg}$ is equal to 0.1 in neutrophil subsystem. As $k_{pg}$ is equal to 0.1, pathogen and binded activated neutrophil diverge at unstable equilibria in neutrophil subsystem. These trends indicate that inflammation oscillation requires the additional intervention or control treatment. Otherwise, the inflammation will constantly occur as time goes by. Fig. 3 (f) and (g) display the phase trajectories in pathogen-free activated neutrophil plane and pathogen-binded activated neutrophil plane, respectively. The arrow in the figure represents the direction of phase trajectory. The stable limit cycles are reach in these phase spaces. Stable limit cycle means that all neighboring trajectories approach the limit cycle as the time approaches infinity [59]. Therefore, Fig. 3 (f) and (g) also reflect that pathogen, free activated neutrophil, binded activated neutrophil will not converge to a stable equilibrium as time approaches infinite when $k_{pg}$ is equal to 0.1. Their values even repeatedly remain in a high level, which will lead to persistent inflammatory. Fig. 3 (h) and (i) display the phase trajectories in (pathogen, free activated neutrophil plane, binded activated neutrophil) space and (pathogen, TNF-$\alpha$, binded activated neutrophil) space, respectively. The stable limit cycles are observed in these two phase space as well.

### 3.2 Bifurcation analysis in monocyte subsystem

Continued bifurcation analysis on the monocyte subsystem are researched. The parameter we analyze in monocyte subsystem is $k_{pg}$. For each parameter, only the system state variables with obvious equilibrium behavior are presented. The bifurcation diagrams of neutrophil subsystem are shown in Fig. 4.

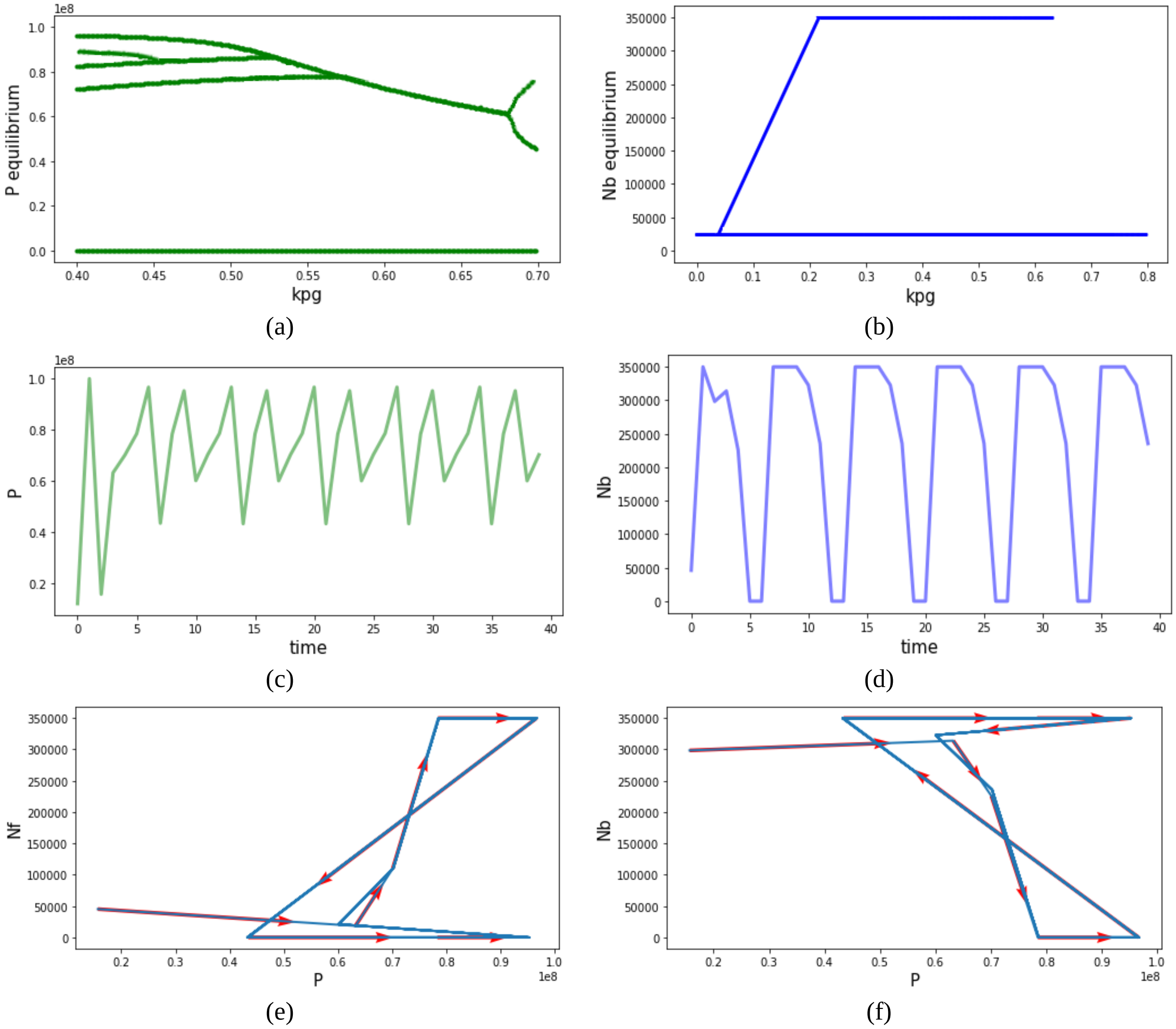

**Fig. 4. (a) Numerical equilibrium curve of pathogen related to parameter $k_{pg}$ in monocyte subsystem.**

**(b) Numerical equilibrium curve of Nb related to parameter $k_{pg}$ in monocyte subsystem.**
**(c) Oscillation behavior of pathogen in monocyte subsystem when $k_{pg}$ is equal to 0.65.**
**(d) Oscillation behavior of Nb in monocyte subsystem when $k_{pg}$ is equal to 0.65.**
**(e) Phase trajectory in P-Nf plane in monocyte subsystem.**
**(f) Phase trajectory in P-Nb plane in monocyte subsystem.**

According to the definition of bifurcation, Fig. 4 (a) and (b) both show the changes in the number of equilibrium and the change in the numerical values of equilibrium when the parameter value is change. In Fig. 4 (a), we observe that more complicate bifurcation behavior of pathogen related to parameter $k_{pg}$ is catch in the monocyte subsystem. There are several bifurcation points representing the change on the number of equilibria. In Fig. 4 (b), when $0 < k_{pg} < 0.1$, there is a bifurcation point leading to the change on the number of $N_b$ equilibrium. Fig. 4 (c) and (d) show the oscillation behaviors of pathogen and $N_b$ when $k_{pg}$ is equal to 0.65. Both pathogen and binded activated neutrophil diverge at unstable equilibria in monocyte subsystem. These trends indicate that in this case the inflammation oscillation will keep happening if there is not any intervention or control measure to change the value of $k_{pg}$. Fig. 4 (e) and (f) display the phase trajectories in pathogen-free activated neutrophil plane and pathogen-binded activated neutrophil plane, respectively. The arrow represents the direction of phase trajectory. The stable limit cycles are observed in these phase spaces. This also can reflect that pathogen, free activated neutrophil and binded activated neutrophil will not converge to a stable equilibrium as time approaches infinite when $k_{pg}$ is equal to 0.65. Their values even repeatedly increase to high level, which will induce the persistent inflammatory.

# 4. Optimal control and RNN-BO optimization algorithm

This section will develop the sepsis model into optimal control model and solve the optimal control strategy using a time-series optimization algorithm named RNN-BO algorithm we detailed proposed in [60]. In this paper, we will consider control strategy variables into sepsis model to represent the level/intensify of sepsis control or intervention treatment strategy. This paper only considers two types of control strategies under two different inflammation situations: one is the control strategy when the load of pathogen remains at high level, but the pro-inflammatory cytokines go down to low level in the early stage of inflammation, and the immune response can't work to the clearance of pathogen; another is the control strategy when the load of pathogen is low, but the immune response is still active.

## 4.1 Control strategy on pathogen and corresponding optimal control model's objective function

Clinical studies show that appropriate antibiotics treatment is effective therapy when the load of pathogen remains at high level, which can effectively control the pathogen replication/growth rate and decrease absolute mortality [11, 12, 13]. Thus, the pathogen growth rate is a controllable parameter. Our sepsis model parameter $k_{pg}$ in Eq. (19) represents the pathogen growth rate (definition provided in Appendix: Table 1). We will consider a $t_f$-dimensional control strategy variable $u_p = \{u_p(t_1), \ldots, u_p(t_f)\}$ to represent the level/intensify of antibiotics treatment control. $u_p(t) \in [u_{pL}, u_{pU}]$ represents the control value at time $t$, $u_{pL}$ and $u_{pU}$ represent the lower bound and upper bound of antibiotics treatment control, respectively. The Eq. (19) will be developed as follows by incorporating the control strategy variable $u_p$:

$$\frac{dP}{dt} = (1 - u_p)k_{pg}P(1 - P^*) - r_{pmk}\frac{[P^n]}{[P^n + k_{c1}^n]}M_{kf}P^* - r_{pn}\frac{[P^n]}{[P^n + k_{c2}^n]}(N_f + N_b)P^* - r_{pm}M_1 - r_{pcd4}\frac{[P^n]}{[P^n + k_{c6}^n]}T_{CD4}P^* - r_{pAb}\frac{[P^n]}{[P^n + k_{c5}^n]}AP^* \quad (39)$$

where $(1 - u_p)$ represents the decrease in pathogen growth rate due to the antibiotics treatment control strategy.

For the optimal control model, the next thing is to determine the objective function for the model. Some good biomarkers/components usually are used as the measure/objective function of immune response or development of sepsis [17]. During the immune response process, the level of pathogen can affect the outcomes of sepsis [61]. M1 macrophage contributes to pathogen clearance but will promote the inflammation as well, M2 macrophage contributes to the removal of apoptotic cells and inhibit inflammation as the same time [27]. M1 macrophage is defined as the up-regulation biomarker and M2 macrophage is the down-regulation biomarker of inflammatory [28]. The ratio of M1/M2 is used as a biomarker correlated with the tissue health status, inflammation associates with higher ratio of M1/M2 [62].

In addition, healthy adaptive immune system plays important role on the recovery of inflammation, CD4+ T cells and CD8+ T cells are two important T cells during adaptive immunity process. CD4+ T cells accelerate the clearance of pathogen [51]. But CD8+ T cells are cytotoxic cells, their primary function is to kill the binded Kupffer cells, binded activated neutrophils, and binded activated monocytes, which will reduce the pathogen clearance ability of immune system [49, 52]. Therefore, the ratio of CD8+ T cell/CD4+ T cell is recognized a biomarker of the ability of adaptive immune system and disease severity, a high CD8+ T cell/CD4+ T cell is associated with increased morbidity and mortality [63, 64].

Therefore, we decide to use the ratio of M1/M2 and the ratio of CD8+ T cell/CD4+ T cell as the objective function when the load of pathogen is high in the early stage of inflammation, to measure the effectiveness of antibiotics treatment control strategy to the development of sepsis. The corresponding objective function is defined as:

$$\min_{u_p \in [u_{pL}, u_{pU}]^{t_f}} w_1 \frac{M1}{M2} + w_2 \frac{T_{CD8}}{T_{CD4}} \tag{40}$$

where $w_1$ and $w_2$ are constant parameters of weight.

### 4.2 Control strategy on TNF-α and corresponding optimal control model's objective function

When the load of pathogen is low in the later stage of inflammation, the inflammation may still present due to the uncontrolled immune response, which will lead to persistent inflammation. If the immune system of host body is weak or uncontrolled, the activated neutrophil cells will still release the toxic chemical substance ROS after finishing pathogen clearance, which is harmful to host tissue and accelerate the death of apoptotic hepatocytes [10]. At the same time neutrophil cells will release TNF-$\alpha$. When TNF-$\alpha$ detects the apoptotic hepatocytes, it will activate and recruit more neutrophil cells to migrate to the site of apoptotic cells. Since phagocytes will constantly attack the host's healthy cells even though there is no pathogen existing in the body, this is a vicious circle to induce persistent infection and eventually develop into server sepsis or organ dysfunction. Some experimental studies show that anti-TNF-$\alpha$ therapy contributes to control the release rate of TNF- $\alpha$ for the above situation, effective anti-TNF-$\alpha$ therapy can improve the outcome of inflammation and save life [15, 16]. No doubt, the release rate of TNF-$\alpha$ is a controllable parameter. Our sepsis model parameter $r_{t2max}$ in Eq. (22) represents the release rate of TNF-$\alpha$ by activated neutrophil (definition provided in Appendix: Table (1). We will consider a $t_f$ -dimensional control strategy variable $u_T = \{u_T(t_1), \ldots, u_T(t_f)\}$ to represent the level/intensify of anti-TNF-$\alpha$ treatment control. $u_T(t) \in [u_{TL}, u_{TU}]$ represents the control value at time $t$, $u_{TL}$ and $u_{TU}$ represent the lower bound and upper bound of anti-TNF-$\alpha$ treatment control, respectively. The Eq. (22) will be developed as follows by incorporating the control strategy variable $u_T$:

$$\frac{dT}{dt} = \left(\frac{r_{t1max} M_{kb}}{m_{t1} + M_{kb}}\right) M_{kb} + \left(\frac{(1-u_T) r_{t2max} N_b}{m_{t2} + N_b}\right) N_b - u_t T \tag{41}$$

where $(1 - u_T)$ represents the decrease in release rate of TNF-$\alpha$ by activated neutrophil due to the anti-TNF-$\alpha$ treatment control strategy.

During the immune response process, the level of inflammatory cytokines both can affect the outcomes of sepsis [61]. TNF- $\alpha$ is a major pro-inflammatory cytokine and IL-10 is a crucial anti-inflammatory cytokine [52]. In many clinical practices, the ratio of TNF- $\alpha$/IL-10 is used as a biomarker to monitor sepsis progression [22, 23, 24]. Therefore, we decide to use the ratio of TNF- $\alpha$/IL-10 as the objective function when the immune response is still active in the later stage of inflammation, to measure the effectiveness of anti-TNF-$\alpha$ treatment control strategy to the development of sepsis. The corresponding objective function is defined as:

$$\min_{u_T \in [u_{TL}, u_{TU}]^{t_f}} \frac{T}{C_A} \tag{42}$$

### 4.3 RNN-BO optimization algorithm

One of our purposes is not only to solve the optimal control that minimizes the objective function value, but also to quickly provide the optimal control strategy when the parameter values or system state variable values are changed. Since the initial system value setting (initial value of sepsis system parameters and system state variables) are associated with the initial conditions of patients, different initial system value settings may lead to different outcomes, and may have different corresponding optimal control strategies. That will waste a lot of time to generate an optimal control strategy if use the optimization algorithm to solve the optimal control model each time for every new given initial values. We know that successful sepsis treatments involve not only optimal dosing of control treatment strategy, the timing of control therapy is important as well [7]. Therefore, an efficient optimization algorithm is key to quickly generate an optimal control strategy, which can reduce the mortality of sepsis and improve patients' quality of care.

The optimization algorithm we use in this paper is named RNN-BO optimization algorithm. The RNN-BO algorithm is a time-series optimization algorithm detailed proposed in our previous paper [60], which combines RNN and an improved BO algorithm. Herein, we briefly introduce the RNN-BO algorithm. The main idea of the RNN-BO algorithm is to use an improved BO algorithm to solve different corresponding low-dimensional optimal control strategies by varying the initial parameter values or system state variable values. Low-dimensional control strategy in here means that the dimension of control strategy what we aim to solve is $d$ $(d < t_f)$ rather than full dimension $t_f$ during this process. The improved BO algorithm is different from the standard BO algorithm. The standard BO algorithm is detailed introduced in [65]. This improved BO algorithm samples the optimal control candidates by combining multi-armed bandit [66] and random search algorithm [67]. Then pick the best solution that minimizes the acquisition function. Acquisition function we used in our RNN-BO algorithm is an approximation function of objective function using lower confidence bound function (LCB) [65]. After the optimization of acquisition function, to increase the solution's accuracy, RNN-BO algorithm does a local search to further optimize this optimal control strategy, which is different from the standard BO algorithm.

For each initial system value setting, we can generate $(t_f - d + 1)$ $d$-dimensional control strategies by using the improved BO algorithm. Since we solve the first $d$-dimensional optimal control strategy start from time 1 to time $d$, the system state variables values over this time period ($t \in [1, d]$) can be calculated based on this first $d$-dimensional optimal control strategy. Then use these system state variables values at time 2 as the initial values, we solve the second $d$-dimensional optimal control strategy start from time 2 to time $d + 1$, the system state variables values over this time period ($t \in [2, d + 1]$) can be calculated based on this second $d$-dimensional optimal control strategy, and so on. If we change the initial parameter value or initial system state variables value, we can generate another $(t_f - d + 1)$ $d$-dimensional control strategies. All of these optimal control strategies are time-series. Store all data pairs (consisting of initial system value settings and corresponding optimal strategies) for further use. For example, if we vary the initial system value setting for $n$ times, then the total number of data pairs we can obtain is $n * (t_f - d + 1)$.

Next, design system value setting as input data and the corresponding optimal control strategy as output. Then use the RNN algorithm to learn those data pairs and generate a model named RNN-BO prediction model. Once provide any initial system value setting, the RNN-BO prediction model can quickly and effectively predict a corresponding $t_f$-dimensional time-series optimal control strategy. The implementation flowchart of the RNN-BO optimization algorithm is shown in Fig. 5.

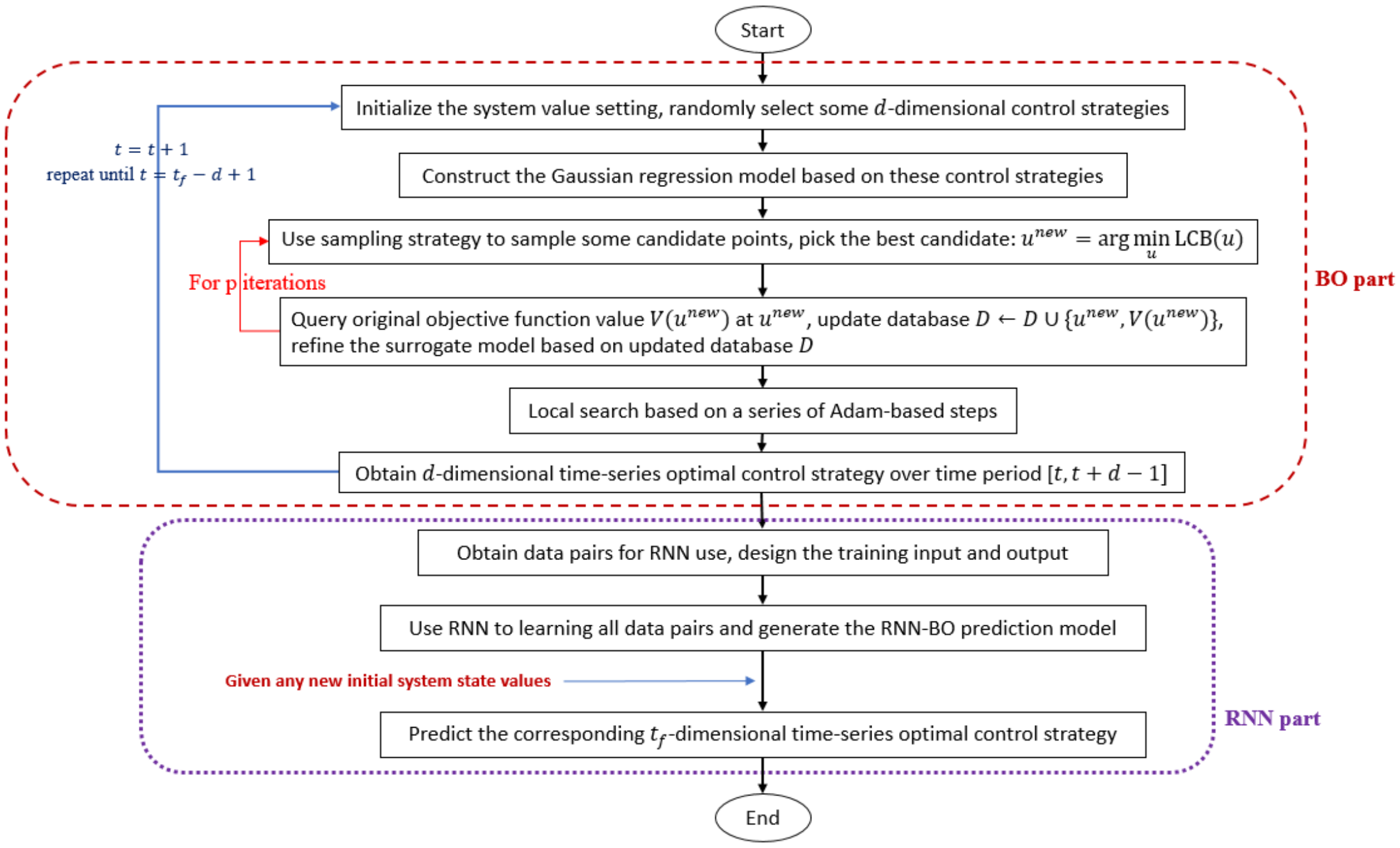

**Fig. 5. Implementation flowchart of the RNN-BO optimization algorithm.**

## 5. Numerical simulation

In this section, we implement numerical simulation tests to solve the optimal control strategy for the sepsis optimal control system in Eqs. (19) – (38) using the RNN-BO algorithm. There are two inflammatory situations that can be controlled as we discussed in Section 4.1 and 4.2: one is when the load of pathogen remains at high level, but the pro-inflammatory cytokines go down to low level in the early stage of inflammation, and the immune response can't work to the clearance of pathogen; two is when the load of pathogen is low, but the immune response is still active. To better demonstrate the effectiveness and efficiency of the RNN-BO algorithm to solve the optimal control strategy on this complex sepsis system, we compare it with the situation without any control, and other two BO algorithm (the standard BO algorithm and a high-dimensional DR-DF BO algorithm proposed in [31])

### 5.1 Numerical results when the optimal control strategy is on pathogen

For the first inflammatory situation, the load of pathogen remains in high level over time, the host's immune system isn't capable to the clearance of pathogen and the pro-inflammatory cytokines go down sooner. In this situation, antibiotics treatment is the effective therapy to control the pathogen replication/growth rate [11, 12, 13]. The first situation can be shown as Fig. 6 (a). TNF-$\alpha$ is an important pro-inflammatory cytokine during immune response process. We can see that when the pathogen goes up in the early stage of inflammation, the immune response is activated. But TNF-$\alpha$ sharply goes down, it means that the macrophages that are responsible to the clearance of pathogen couldn't be recruited and activated. In this case, the load of pathogen will remain in high level, this may lead to the death due to pathogen infection.

When the control strategy is antibiotics treatment control strategy, the objective function is to minimize the sum of ratio of M1/M2 and the ratio of CD8+ T cell/CD4+ T cell. Higher ratio is associated with severe inflammation. The simulation results are shown in Fig. 6. The running time of the standard BO algorithm to generate the optimal control strategy is about 45 seconds. The running time of the DR-DF BO algorithm is about 25 seconds. But the RNN-BO algorithm is different from other two algorithms, it learns the historical data. Once the RNN-BO prediction model is ready, the RNN-BO algorithm only takes about 2 seconds to predict the optimal control strategy by giving the same initial system values as the other two BO algorithms.

Fig. 6 (b) shows the control strategies from three algorithms. The optimal control strategy of standard BO algorithm performs obvious fluctuation over time. The optimal control strategy of DR-DF BO algorithm is more stable. The optimal control predicted from RNN-BO algorithm is lower at the early stage of inflammation, then become high level at the later stage of inflammation when it recognizes the load of pathogen is still in high level. According to the trends of the optimal control strategies, the optimal control strategy predicted by RNN-BO algorithm may be more reasonable.

Fig. 6 (c) shows the comparison on ratio of CD8+ T cell/CD4+ T cell over time. We can see that the ratio when the system is without control is significantly higher than the ration when the system is with control. From the smaller figure in Fig. 6 (c), after applying the optimal control strategies generated by the standard BO algorithm, DR-DF BO algorithm, and RNN-BO algorithm, the ratios of CD8+ T cell/CD4+ T cell perform the same trends. That means those three algorithms reach similar optimization performances on this ratio for our sepsis optimal control system, they all have effective impact on controlling the inflammation. Fig. 6 (c) shows the comparison on ratio of M1/M2 over time. We can see that the ratio when the system is without control is also significantly higher than the ration when the system is with control. The ratio without control will gradually increase up to 40,000. From the smaller figure in Fig. 6 (d), the ratios with control are effectively controlled at low level. The DR-DF BO algorithm and RNN-BO algorithm have similar great performance, both slightly outperform the standard BO algorithm on the ratio of M1/M2. Fig. 6 (e) shows the accumulated objective function values over time of different methods. Since the objective function is the sum of ratio of CD8+ T cell/CD4+ T cell and ratio of M1/M2, the trends of accumulated objective function value are like the trends of the ratios.

According to Fig. 6, taking antibiotics treatment control is necessary to control the progression of inflammation when the load of pathogen is high in the early stage of inflammation. Overall, the optimal control predicted by the RNN-BO algorithm is slightly better than the standard BO algorithm and DR-DF BO algorithm with only 2 seconds running time.

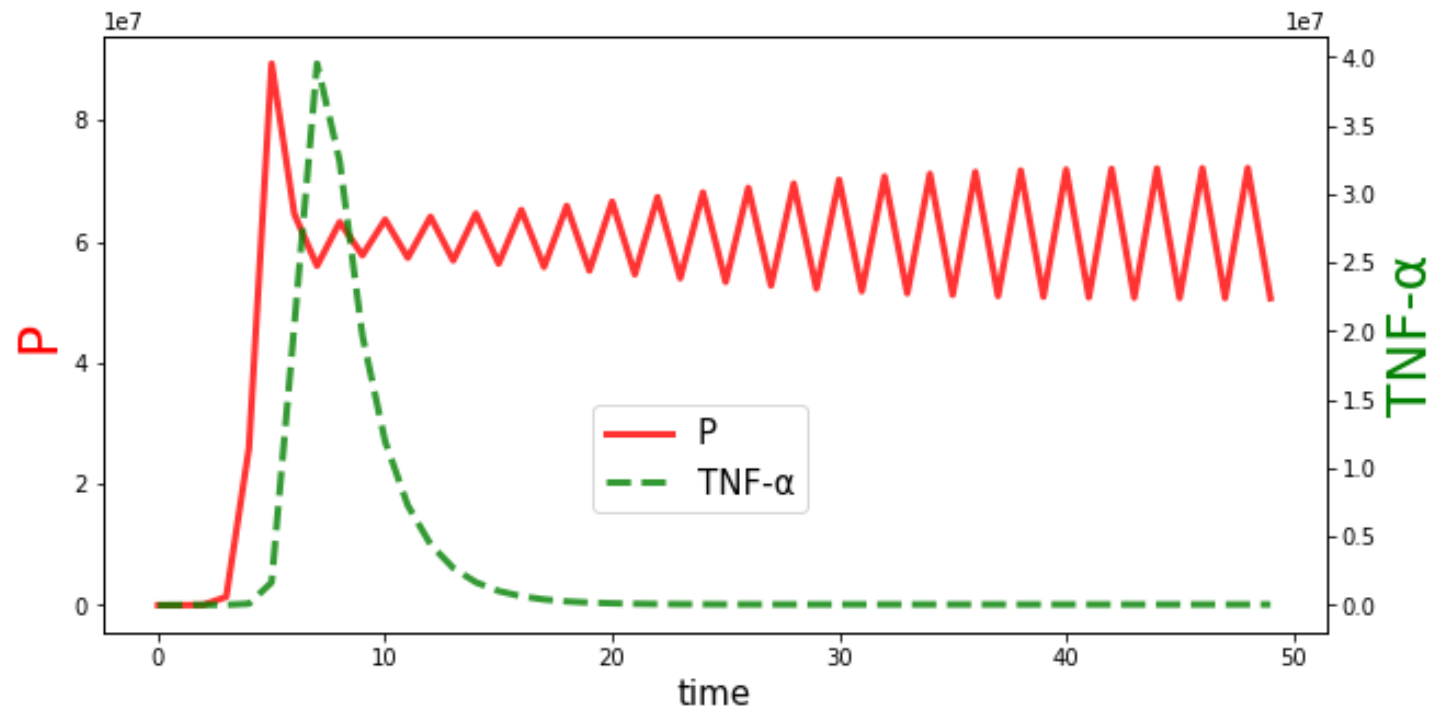

(a)

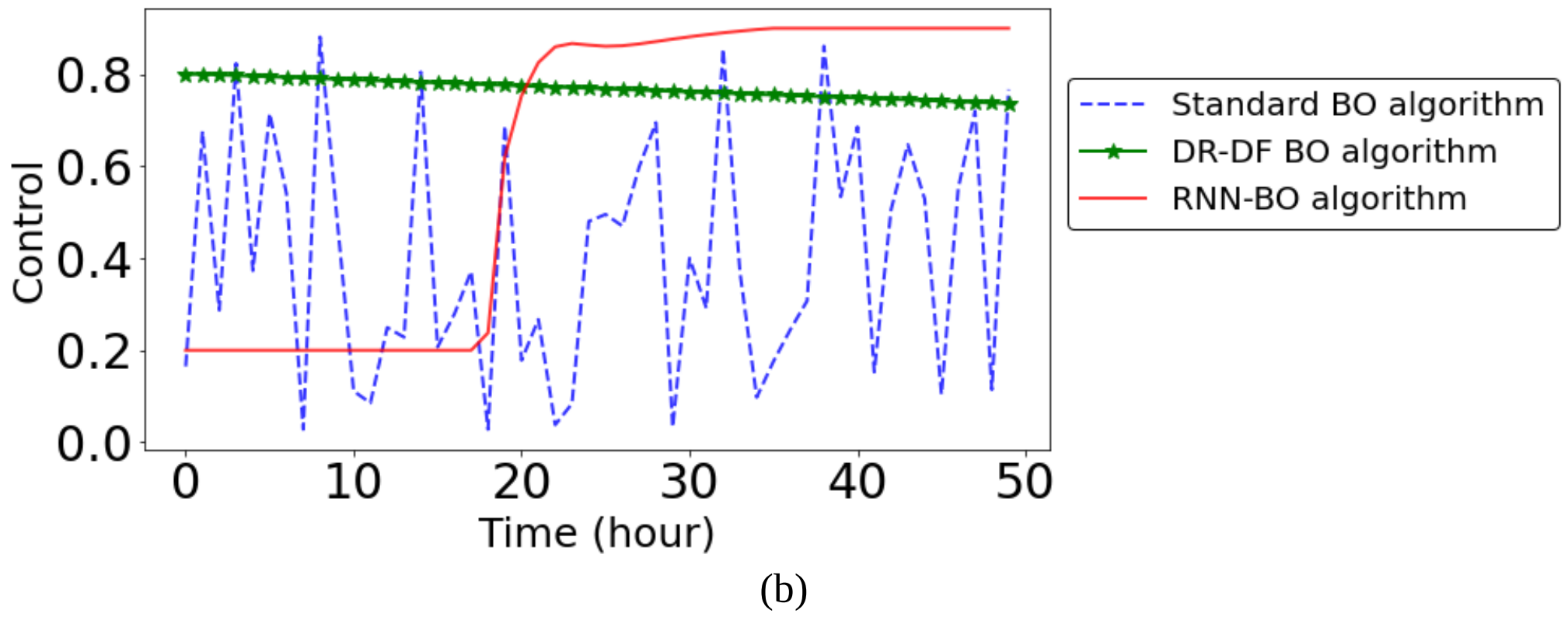

(b)

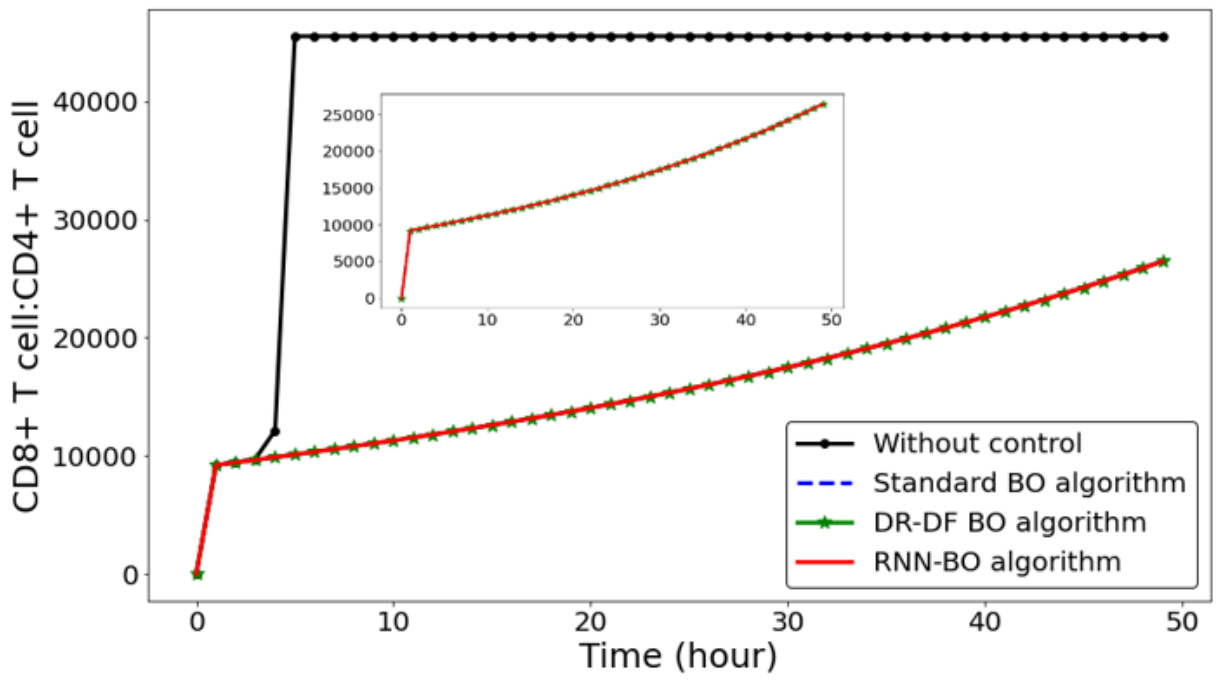

(c)

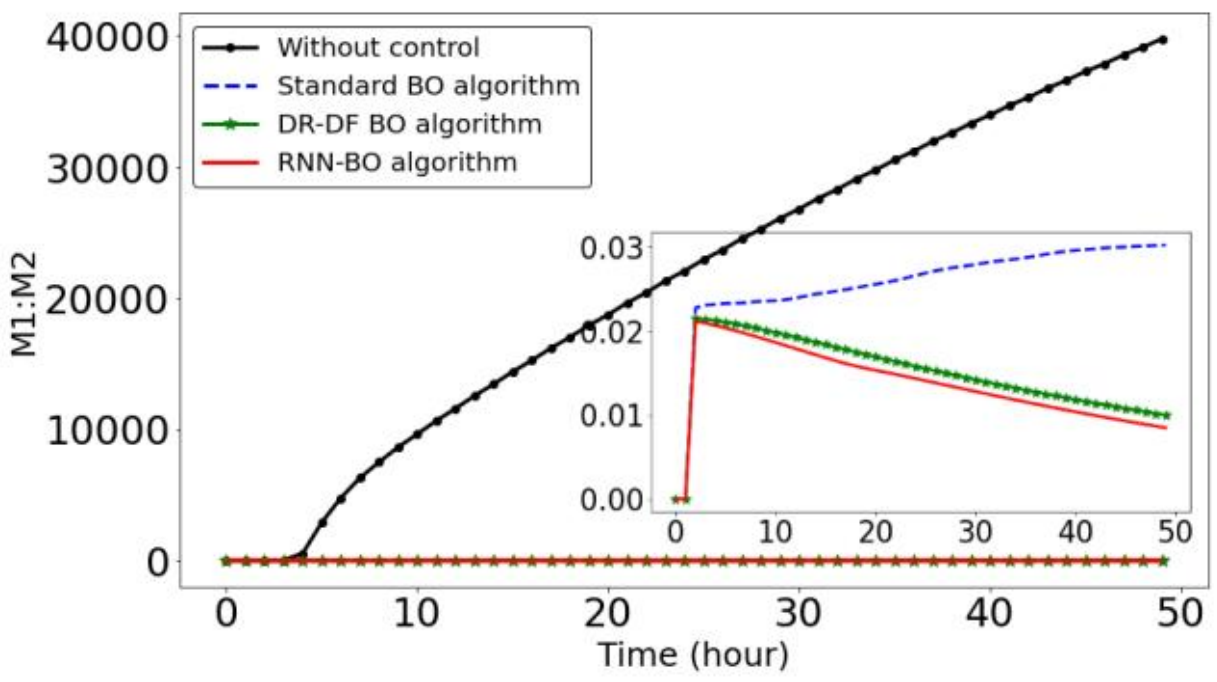

(d)

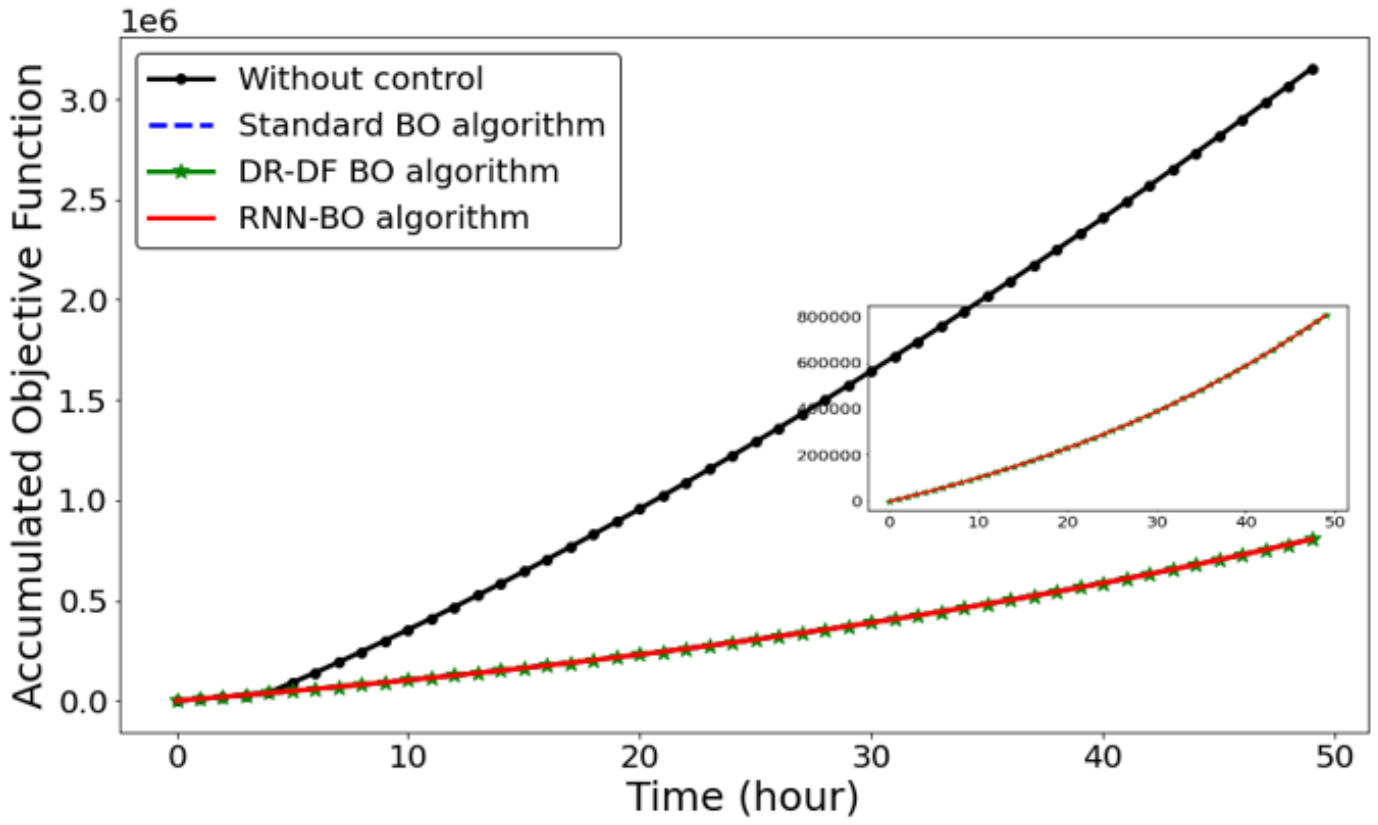

(e)

**Fig. 6. (a) Trends of Pathogen and TNF- $\alpha$ in the first inflammatory situation.**
**(b) Optimal control strategies of different optimization algorithms in the first inflammatory situation.**
**(c) Ratio of $\frac{M1}{M2}$ comparison by different algorithms in the first inflammatory situation**.
**(d) Ratio of $\frac{T_{CD8}}{T_{CD4}}$ comparison by different algorithms in the first inflammatory situation**.
**(e) Accumulated objective function over time of different algorithms in the first inflammatory situation.**

## 5.2 Numerical results when the optimal control strategy is on TNF-α

For the second inflammatory situation, the immune response is still active when the load of pathogen is low. This means that macrophages constantly attack the host's healthy cells after they finish the clearance of pathogen. In this situation, the host will perform persistent inflammation and tend to develop into organ dysfunction. From previous clinic practices, anti-TNF-$\alpha$ therapy is an effective control treatment to the second situation [15, 16]. The second situation can be shown as Fig. 7 (a). We can see that the load of pathogen grows up quickly in the early stage of inflammation, the immune response is activated. After the load of pro-inflammatory cytokine TNF- $\alpha$ increases, the pathogen starts to go down until all pathogens are cleared. However, after all pathogens are cleared, the load of TNF- $\alpha$ remains in a high level. The immune response keeps active even there is no pathogen in the host's body. This case may lead to the death due to persistent inflammation.

When the control strategy is anti-TNF-$\alpha$ treatment control strategy, the objective function is to minimize the ratio of TNF- $\alpha$/IL-10. The simulation result is shown in Fig. 7. Higher ratio is associated with severe inflammation. The running times of different algorithms in this section are similar as they performed in Section 5.1.

Fig. 7 (b) shows the control strategies from three algorithms. The optimal control strategy of standard BO algorithm performs obvious fluctuation over time. The optimal control strategy of DR-DF BO algorithm is more stable. The optimal control predicted from RNN-BO algorithm is lower at the early stage of inflammation, then sharply increase to a high level when it recognizes the load of TNF- $\alpha$ keeps increasing even the pathogen has already started to decrease. According to the trends of the optimal control strategies, the optimal control strategy predicted by RNN-BO algorithm may be more reasonable.

Fig. 7 (c) shows the comparison on ratio of TNF- $\alpha$/IL-10 over time. We can see that the ratio when the system is without control is significantly higher than the ration when the system is with control. The ratio without control will gradually increase up to $2.5 \times 10^9$. From the smaller figure in Fig. 7 (c), the ratios with control are effectively controlled, the highest ratios with control are about $10^3$ times lower than the highest ratio of without control. The RNN-BO algorithm outperforms the standard BO algorithm and the DR-DF BO algorithm. Fig. 7 (d) shows the accumulated objective function values over time of different methods. Since the objective function in the second inflammatory situation is the ratio of M1/M2, the trend of accumulated objective function value performs like the trends of the ratio.

According to Fig. 7, taking anti-TNF-$\alpha$ treatment control is necessary to control the progression of inflammation when the load of pathogen is low but the immune response is still active in the later stage of inflammation. Overall, the optimal control generated by the RNN-BO algorithm is better than the standard BO algorithm and DR-DF BO algorithm.

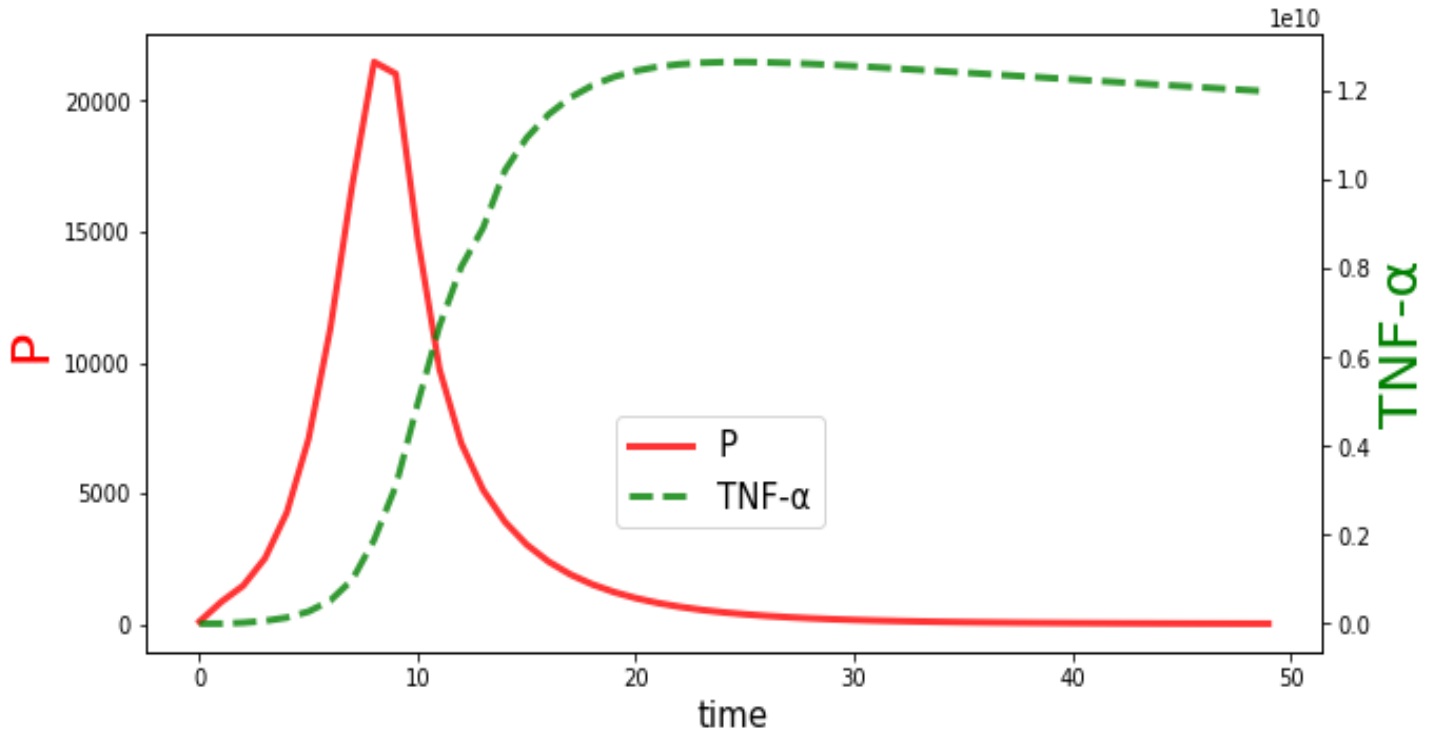

(a)

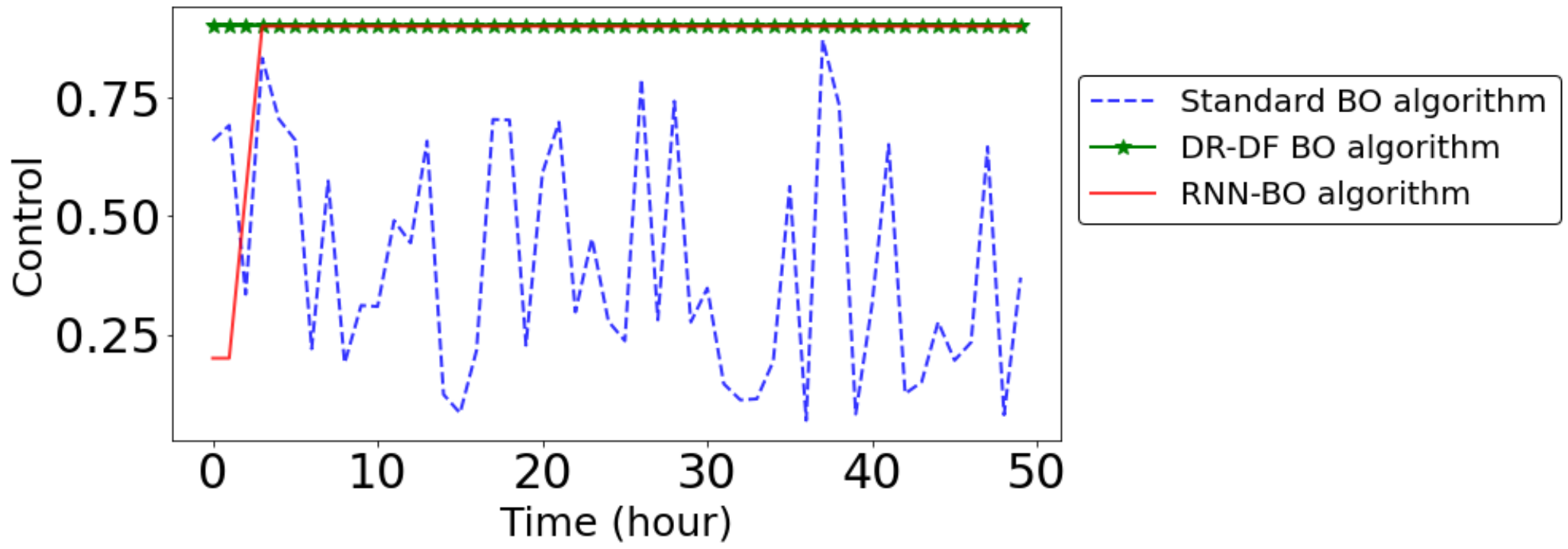

(b)

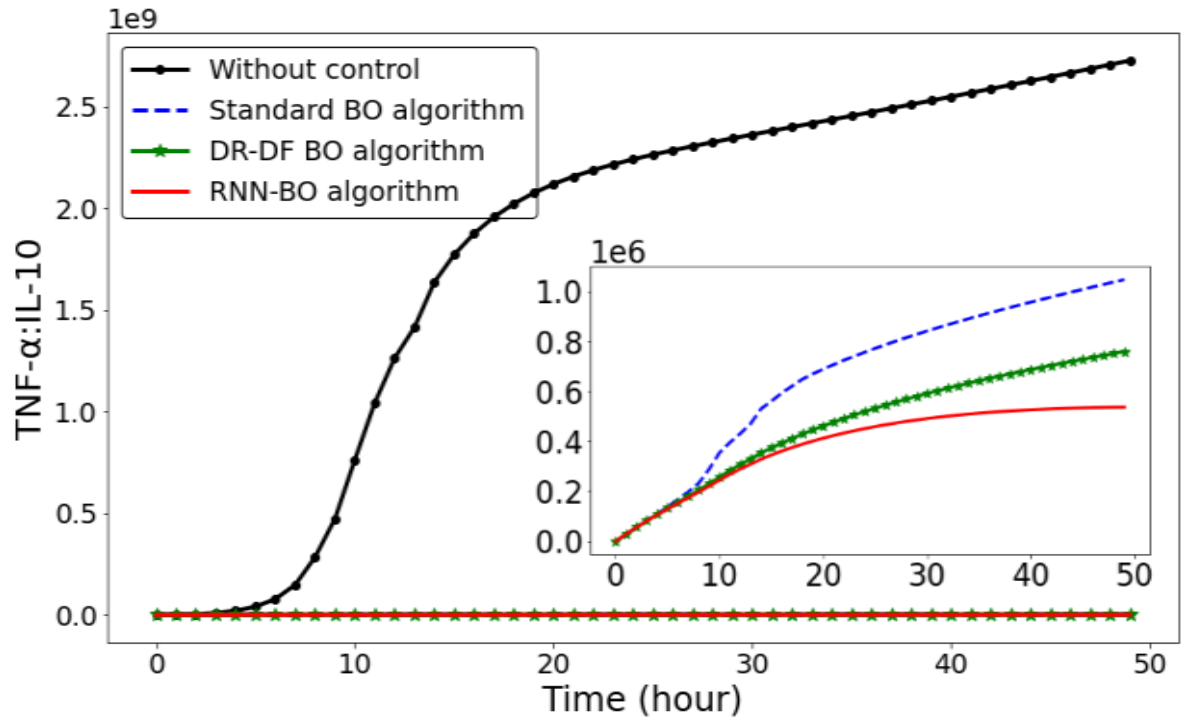

(c)

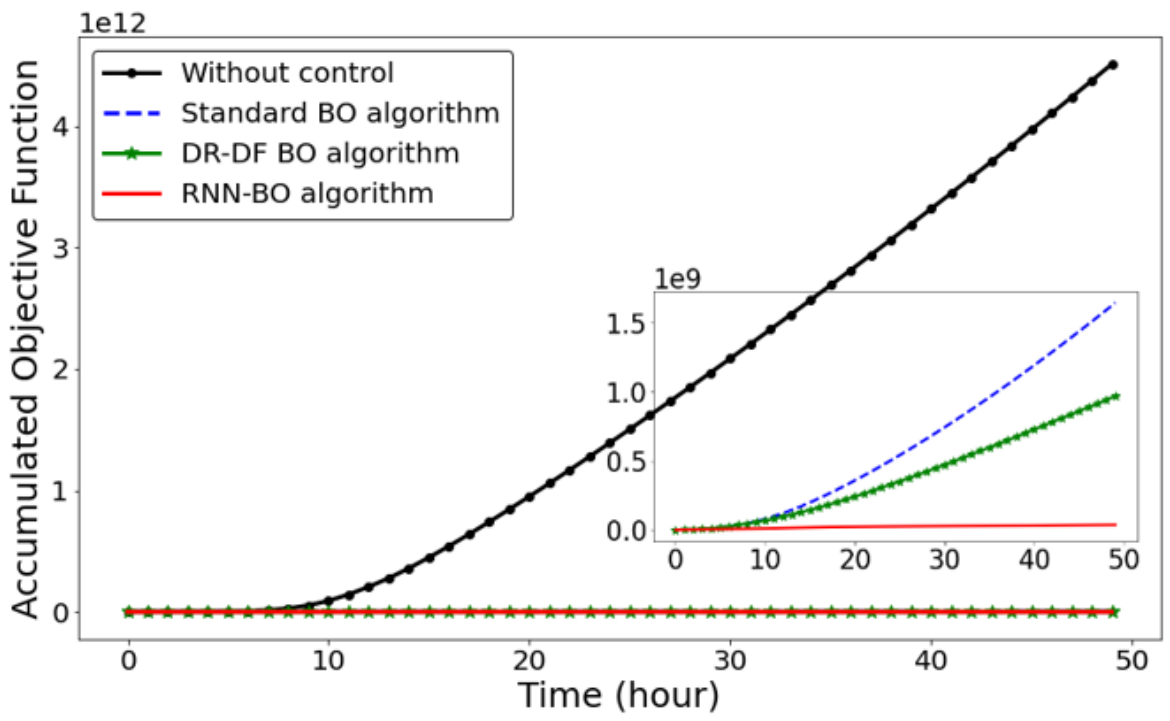

(d)

**Fig. 7. (a) Trends of Pathogen and TNF- $\alpha$ in the second inflammatory situation.**
**(b) Optimal control strategies of different optimization algorithms in the second inflammatory situation.**
**(c) Ratio of $\frac{T}{C_A}$ comparison by different algorithms in the second inflammatory situation**.
**(d) Accumulated objective function over time of different algorithms in the second inflammatory situation.**

## 6. Conclusion and future work

This paper improves a complex nonlinear sepsis model on the monocyte part and adaptive immunity part, which is better to study the progression of the delicate immune response system. The bifurcation analysis of our sepsis subsystem presents the model behaviors under some system parameters, but also shows the necessary of control treatment and intervention therapy for the sepsis development. If the sepsis system is without considering any control treatment under some parameter and initial system value settings, the system will perform persistent inflammation outcomes (harmful infection oscillation outcomes) as time goes by. Thus, this paper develops the improved nonlinear sepsis model into an optimal control system. According to some existing clinic practices, this paper determines to apply some authorized and recommended sepsis biomarkers as our objective function of studied sepsis optimal control system to measure the development of sepsis. Next, an RNN-BO optimization algorithm is introduced to predict the optimal control strategy. The most advantage of RNN-BO algorithm is that it learns the historical optimal control strategy and generates a prediction model. Once there is a new sepsis patient with different initial condition (is associated with initial system value setting), the RNN-BO algorithm is capable to predict the corresponding optimal control strategy for this patient in short time. Some comparison simulation experiments with other optimization algorithms are carried out. Simulation results demonstrate the effectiveness and efficiency of the RNN-BO algorithm on driving the optimal control solution for a complex nonlinear sepsis optimal control system. As the healthcare field develops, the mathematical study and optimal control research of sepsis will continue to grow. To better express sepsis via mathematical model, and propose more effective optimization algorithm for providing the optimal control strategy to improve quality of clinic therapy or reduce the mortality of sepsis, are both our further research directions.

# Appendix

**Table 1. Definition and experimental simulation values of parameters**

| Parameter | Definition | Value | Reference |
|---|---|---|---|
| $k_{pg}$ | Pathogen growth rate | 0-3.6/h | [68] |
| $P_\infty$ | Pathogen carrying capacity | $10^8$ cells | [69] |
| $r_{pmk}$ | Rate at which pathogens are killed by Kupffer cells | 0.03/per kupffer cell/h | [70] |
| $n$ | The extent of pathogen binding to Kupffer cells | 2 | [10] |
| $k_{c1}$ | Number of Kupffer cells which phagocytose half of pathogen | 0.03 cells/h | [70] |
| $r_{pn}$ | Rate at which pathogens are killed by neutrophils | 20-100/per neutrophil/h | [71] |
| $k_{c2}$ | Concentration of neutrophils which phagocytose half of pathogen | $1.5 \times 10^{-4}$/h | [72] |
| $k_{mk}$ | Proliferation rate of Kupffer cells under inflammation | 0.015 – 2/h | [10] |
| $K_\infty$ | Kupffer cells carrying capacity | (16-20)$\times 10^6$ cells/g liver | [73] |
| $k_{mkub}$ | Unbinding rate of binding Kupffer cells | 0.1-0.77/h | [74] |
| $u_{mk}$ | Killing rate of free Kupffer cells induced by binding to pathogen | 0.23-0.9/h | [74] |
| $r_{t1max}$ | The maximum number of TNF- $\alpha$ being released by Kupffer cells per enzyme molecule per hour | 10/h | [10] |
| $m_{t1}$ | Number of Kupffer cells at which the reaction rate is half of maximal production rate | 10000 cells | [10] |
| $r_{t2max}$ | The maximum number of TNF- $\alpha$ being released by neutrophils per enzyme molecule per hour | 1000/h | [10] |
| $m_{t2}$ | Number of activated neutrophils at which the reaction rate is half of maximal production rate | 10000 cells | [10] |
| $u_t$ | Degradation rate of TNF- $\alpha$ | 0.025-0.5/h | [75] |
| $k_{rd}$ | Influx rate of neutrophils into blood vessel | 0.1-0.72/h | [76] |
| $N_S$ | Maximum amount of neutrophils in liver | $3.5 \times 10^5$/h | [77] |
| $u_{nr}$ | Apoptotic rate of resting neutrophils | 0.069-0.12/h | [78] |
| $k_{nub}$ | Unbinding rate of activated neutrophils | 0.01-0.5/h | [10] |
| $u_n$ | Apoptotic rate of activated neutrophils | 0.05/h | [78] |
| $k_{r1}$ | Auxiliary parameter associated with the activation rate of resting neutrophils | 3/h | [10] |
| $u_{r1}$ | Degradation rate of parameter $r_1$ to maintain a slow-saturation curve | 0.003/h | [10] |
| $r_{hn}$ | Rate at which activated neutrophils kill apoptotic hepatocytes | 9000/per neutrophil/h | [10] |
| $k_{c3}$ | Concentration of activated neutrophils which phagocytose half of apoptotic hepatocytes | 0.04 cells/h | [10] |
| $A_\infty$ | Number of hepatocytes in liver | $3.2 \times 10^8$ cells/h | [10] |
| $r_{ah}$ | Recovery rate of apoptotic hepatocytes | 0.5-2/h | [79] |
| $C_\infty$ | Dissociation rate of IL-10 | 0.02 | [10] |
| $u_{mn}$ | Rate at which activated neutrophils are killed by inflammatory monocytes | 200/monocyte/h | [10] |
| $k_{mr}$ | Influx rate of monocytes into blood vessel | 0.5/h | [80] |
| $M_S$ | Resting monocyte carrying capacity in blood vessel | 50000 cells | [81] |
| $r_2$ | Influx rate of monocytes in liver | 80/h | [82] |
| $u_{mr}$ | Apoptotic rate of resting monocytes | 0.2 | [10] |
| $u_m$ | Apoptotic rate of activated monocytes | 0.08 | [83] |
| $k_{umb}$ | Unbinding rate of binding activated monocytes | 0.4 | [84] |
| $r_{pm}$ | Rate at which pathogens are killed by inflammatory monocytes | 7/monocyte/h | [85] |
| $k_{c4}$ | Number of monocytes that phagocytose half of pathogen | 0.002 cells/h | [85] |
| $r_{h1max}$ | The maximum number of HMGB-1 being released by monocytes per enzyme molecule per hour | 0.001 | [10] |
| $mh_1$ | Number of monocytes generate half of maximal HMGB-1 production rate | 10,000 | [10] |
| $u_h$ | Degradation rate of HMGB-1 | 0.5–3 | [10] |
| $r_{camax}$ | The maximum number of IL-10 being released by monocytes per enzyme molecule per hour | 10,000 | [10] |
| $C_{Ah}$ | Number of monocytes generate half of maximal HMGB-1 production rate | 10,000 | [10] |
| $u_{ca}$ | Degradation rate of IL-10 | 0.02 | [10] |
| $r_{pcd4}$ | Rate at which pathogens are killed by CD4+ T cells | 8 | [85] |
| $k_{c5}$ | Concentration of antibody which kills half of pathogen | 0.035 | [10] |
| $k_{c6}$ | Concentration of CD4+ T cells which kill half of pathogen | 0.0015 | [10] |
| $r_{pAb}$ | Rate at which pathogens are killed by antibody | 1 | [10] |
| $r_{Mkbcd8}$ | Rate at which binding Kupffer Cells are killed by CD8+ T cells | 0.25 | [86] |
| $r_{Nbcd8}$ | Rate at which binding activated neutrophils are killed by CD8+ T cells | 0.25 | [86] |
| $k_{c7}$ | Concentration of CD8+ T cells which kill half of binding antigen presenting cell | 0.0015 | [10] |
| $r_{cd4Mb}$ | Rate at which CD4+ T cells bind to activated monocytes | 4 | [86] |
| $r_{cd8Mb}$ | Rate at which CD8+ T cells bind to activated monocytes | 4 | [86] |
| $k_{c8}$ | Activated monocyte concentration produces half occupation on T cells | 0.0075 | [10] |
| $r_{Mbcd8}$ | Rate at which binding activated monocytes are killed by CD8+ T cells | 0.25 | [86] |
| $k_{cd4M}$ | Rate at which binding CD4+ T cells are killed by activated monocytes | 0.73–2 | [87] |
| $k_{cd8M}$ | Rate at which binding CD8+ T cells are killed by activated monocytes | 0.73–2 | [87] |
| $k_{c9}$ | B cell concentration produces half occupation on T cells | 0.045 | [10] |
| $k_{c10}$ | Concentration of activated monocytes which kill half of binding T cells | 0.018 | [10] |
| $k_{cd4}$ | The influx rate of CD4+ T cells to blood vessel | 0.014 | [88] |
| $T_{CD4\infty}$ | CD4+ T cell carrying capacity in the blood vessel | $27.4 \times 10^6$ | [88] |
| $u_{cd4}$ | Degradation rate of CD4+ T cells | 0.00083–0.001 | [88] |
| $k_{cd8}$ | The influx rate of CD8+ T cells to blood vessel | 0.0625 | [88] |
| $T_{CD8\infty}$ | CD8+ T cell carrying capacity in the blood vessel | $5 \times 10^6$ | [88] |
| $u_{cd8}$ | Degradation rate of CD8+ T cells | 0.00079–0.001 | [88] |
| $k_B$ | The influx rate of B cells to blood vessel | 0.0122 | [88] |
| $B_\infty$ | B cell carrying capacity in the blood vessel | $28.6 \times 10^6$ | [88] |
| $r_{Bt}$ | Rate at which B cells bind to T cells | 1–10 | [10] |
| $u_B$ | Degradation rate of B cells | 0.00012–0.00016 | [89, 90] |
| $r_{Abmax}$ | The maximum production amount of antibody by B cells | 0.00053 | |
| $m_{Ab}$ | Number of B cells at which the reaction rate is half of maximum production rate | 10,000 | [10] |
| $u_{Ab}$ | Degradation rate of antibody | 0.0035–0.01 | [91] |